\documentclass[twoside]{article}
\usepackage[usenames,dvipsnames]{color}
\usepackage{graphicx}
\usepackage{camnum}
\usepackage{harvard}

\usepackage{amsmath}
\usepackage{eufrak}
    \usepackage{ulem}
    \usepackage{chngcntr}
\usepackage{mathtools}
\usepackage{bm}

\newtheorem{remark}{Remark}
\counterwithin{table}{section}

\usepackage{tikz,pgfplots}
        \pgfplotsset{compat = 1.3}
        \pgfplotsset{minor grid style={dotted}} \pgfplotsset{major grid
        style={dashed}}

        \pgfplotsset{every x tick label/.append style={font=\footnotesize,
        yshift=0.25ex}}

        \pgfplotsset{every y tick label/.append
        style={font=\footnotesize, xshift=0.25ex}}

        \definecolor{colorclassyorange}{rgb}{0.95000,0.32500,0.09800}
        \definecolor{colorchromeyellow}{rgb}{1.00000,0.6549,0}%
        \definecolor{colorpaleyellow}{rgb}{1.00000,0.8549,0.1}%

        \definecolor{colorclassyblue}{rgb}{0.00000,0.44706,0.74118}%
        \definecolor{colorpurple}{rgb}{0.49400,0.18400,0.55600}%
        \definecolor{colorfuschia}{rgb}{0.95039,0.0,0.95039}%
        \definecolor{colorlemongreen}{rgb}{0.6,0.8,0}%
        \definecolor{colorreal}{rgb}{0.92941,0.79412,0.12549}%
        \definecolor{colorimag}{rgb}{0.00000,0.49804,0.00000}%
        \definecolor{colorabs}{rgb}{1.00000,0.00000,0.00000}%

\newcommand{\ee}{{\mathrm e}}
\newcommand{\ii}{{\mathrm i}}

\newcommand{\sBCH}{\CC{sBCH}}
\newcommand{\K}{\mathcal{K}}

\newcommand{\FFF}{\mathcal{F}}

\newcommand{\DDD}{\mathcal{D}}
\newcommand{\WW}{\mathcal{W}}
\newcommand{\dx}{\partial_x}

\def\O#1{{\cal O}\left(#1\right)}

\newcommand{\Int}[4]{\int_{#2}^{#3}\!#4\,\mathrm{d}#1}
\newcommand{\Ang}[2]{\left\langle #2 \right\rangle_{#1}}

\newcommand{\ang}[2]{\left\langle #1 \right\rangle_{#2}}

\newcommand{\ve}{\varepsilon}
\newcommand{\schr}{Schr\"odinger }

\newcommand{\Hc}[1]{#1}

\newcommand{\dt}{\Hc{h}}

\newcommand{\mTh}[1]{\Hc{\Theta_{#1}(\dt)}}

\newcommand{\mom}[2]{\Hc{\mu_{#1,#2}(\dt)}}
\newcommand{\momO}{\Hc{\mu_{0,0}(\dt)}}
\newcommand{\Lf}[3]{\Hc{\Lambda\left[#1\right]_{#2,#3}\!(\dt)}}

\newcommand{\norm}[2]{\left\| #2 \right\|_{#1}}

\newcommand{\momS}[3]{\Hc{\mu_{#1,#2,#3}(\dt)}}
\newcommand{\momOS}[1]{\Hc{\mu_{#1,0,0}(\dt)}}
\newcommand{\LfS}[4]{\Hc{\Lambda\left[#1\right]_{#2,#3,#4}\!(\dt)}}

\numberwithin{equation}{section}

\begin{document}

\title{Solving \schr equation in semiclassical regime with
highly oscillatory time-dependent potentials}

\author{Arieh Iserles,\footnote{Department of Applied Mathematics and Theoretical Physics,
 University of Cambridge, Wilberforce Rd, Cambridge CB3 0WA, UK.}\ \
  Karolina Kropielnicka\footnote{Institute of Mathematics, Polish Academy of Sciences,
  8 \'Sniadeckich Street, 00-656 Warsaw, Poland.}\ \  \&
  Pranav Singh\footnote{Mathematical Institute, Andrew Wiles Building, University of Oxford, Radcliffe Observatory Quarter, Woodstock Rd, Oxford OX2 6GG, UK and Trinity College, University of Oxford, Broad Street, Oxford OX1 3BH, UK.}}
\maketitle

\vspace*{6pt}

\vspace*{6pt} \noindent   \textbf{Keywords:} Schr\"odinger equation, time
dependent potentials, semiclassical regime, highly oscillatory potentials,
large time steps, integral preservation, simplified commutators, Magnus
expansion, symmetric Zassenhaus splittings, Lanczos iterations




\begin{abstract}
\schr equations with time-dependent potentials are of central importance in
quantum physics and theoretical chemistry, where they aid in the simulation
and design of systems and processes at atomic and molecular scales.
Numerical approximation of these equations is particularly difficult in the
semiclassical regime because of the highly oscillatory nature of solution.
Highly oscillatory potentials such as lasers compound these difficulties
even further. Altogether, these effects render a large number of standard
numerical methods less effective in this setting. In this paper we will
develop a class of high-order exponential splitting schemes that are able
to overcome these challenges by combining the advantages of
integral-preserving simplified-commutator Magnus expansions with those of
 symmetric Zassenhaus splittings. This allows us to use large time steps
in our schemes even in the presence of highly oscillatory potentials and
solutions.
%

  %
%
%
\end{abstract}

\newpage

\setcounter{section}{0}

\section{Introduction}
The \schr equation in the semiclassical regime with a time-dependent
electric field is the initial value problem,
\begin{equation}\label{eq:schr}
\partial_t u(x,t)=\ii \ve \dx^2 u(x,t) - \ii \ve^{-1} V(x,t) u(x,t),\quad x\in I \subseteq \BB{R},\ t \geq 0,\ u(x,0)=u_0(x).
\end{equation}
Here $V$ is a real-valued, time-dependent electric field and $0<\ve \ll 1$ is
the {\em semiclassical} parameter whose small size induces rapid oscillations
of solution which cause obvious difficulties with numerical discretisation.
The regularity that we require in the potential and the solution varies with
the accuracy of a specific method (in the family of high-order methods that
we will introduce here). For the sake of simplicity, however, we assume that
they are sufficiently smooth -- $V(\cdot,t) \in \CC{C}_p^\infty(I; \BB{R})$
and $u(\cdot,t) \in \CC{C}_p^\infty(I; \BB{C})$.

Effective numerical approximation of this fundamental equation of quantum
mechanics is of great importance and presents many computational challenges.
Among the more challenging applications, for instance, is optimal control of
a quantum system whose dynamics are governed by \R{eq:schr}. This is the
inverse problem of determining a time-dependent electric field or laser
$V(x,t)$ that
 steers the initial system of underlying particles
to a desired state. In practice, this involves solving a variational problem
which, in turn, requires repeated numerical approximations of forward and
backward equations of the form \R{eq:schr} in conjunction with an
optimisation algorithm.

Such applications naturally require highly accurate but computationally
inexpensive numerical methods. A wide range of highly effective methods has,
indeed, been developed for solving the \schr equation with time-dependent
potentials,  a very active research problem in theoretical chemistry, quantum
physics and numerical mathematics
\cite{TalEzer1992,Peskin,SanzSerna96,TremblayCarrington04,NdongEzerKosloff10,alvermann2011hocm,BlanesCasasMurua17,blanes17quasimagnus,Schaefer,IKS18sinum}.

In the setting under consideration, however, the complexity of the task is
compounded due to two factors. First of all, the oscillatory nature of the
solution in quantum mechanics is significantly more pronounced in the
semiclassical regime. Secondly, these applications often involve potentials
that feature high temporal oscillations (such as lasers).

Lanczos-based methods, e.g.\ \cite{IKS18sinum}, are very effective in the
atomic regime characterised by $\ve=1$, less so in the
semiclassical regime where $\ve \ll 1$. Although they can still be
applied in this context, they become computationally expensive. This happens
because the superlinear accuracy, the most welcome feature of Krylov subspace methods such as Lanczos
iterations, commences only once the number of iterations has
exceeded the spectral radius of the exponent,\footnote{To be more precise, the
number of Lanczos iterations needs to exceed half the length of the interval
over which the eigenvalues of the skew-Hermitian exponent are spread
\cite{hochbruck97ksa}.}  which grows like $\O{h \ve^{-1}}$ in the
semiclassical regime. Effectively, this forces us to use either a large
number of Lanczos iterations or a severe depression of the time step, $h$,  to attain sufficient accuracy. In either case, the cost grows
prohibitively.


%
%

This growth of the exponent is due to the $\O{\ve}$ wavelength
spatio-temporal oscillations that appear in the solution of this equation
regardless of the smoothness of the initial conditions
\cite{bao02ots,jin11mcm,Singh17}. In order to correctly resolve these
oscillations, we require
\begin{equation}
\label{eq:Nve} N = \O{\ve^{-1}}
\end{equation}
grid points in the spatial discretisation \cite{bao02ots}, causing the
discretisation\footnote{Throughout this manuscript we use the notation
$\mathcal{L} \leadsto L$ to denote spatial discretisation of an operator
$\mathcal{L}$ by the matrix $L$.} of $\ii \ve \dx^2$ to scale as
$\O{\ve^{-1}}$ since
\begin{equation}
\label{eq:dxve}
\dx^k \leadsto \K_k = \O{N^k} =
\O{\ve^{-k}}. \end{equation}


Recently a specialized approach for semiclassical regime was developed in
\cite{bader14eaf} and \cite{BIKS}. In \cite{bader14eaf} a new methodology of
{\em symmetric Zassenhauss splitting} was introduced for effectively treating
the case of time-independent potentials. This splitting was termed an {\em
asymptotic} splitting due to the manner in which it separates vastly
differing scales, enabling far cheaper computation of exponentials in the
semiclassical regime. In \cite{BIKS} we extended this approach to the case of
time-dependent potentials in semiclassical regime by applying the symmetric
Zassenhaus splitting algorithm (see Appendix~\ref{sec:Zassenhaus}) to the
Magnus expansion \cite{blanes09tme,iserles00lgm} and using ideas pioneered in \cite{mko99fla}. Since this approach commences from a
Magnus expansion where the integrals are discretised at the outset, however,
it can be ineffective in the case of highly oscillatory potentials.

%

\subsection{The  approach of this paper}
In this paper we develop an effective computational approach for \schr equations in
semiclassical regime featuring time-dependent potentials that may be highly
oscillatory. We choose optimal tools for this purpose, introducing new ones
where required for further  efficiency. In particular:

\begin{enumerate}
\item We start the derivation of our schemes from the integral-preserving
    simplified-commutator Magnus expansions of \cite{IKS18sinum}.
\item Instead of directly resorting to Lanczos iterations for
    exponentiating the Magnus expansion (which proves to be a costly
    approach in the semiclassical regime), we perform a symmetric
    Zassenhaus splitting on this Magnus expansion. The use of this asymptotic
    splitting algorithm separates scales. Consequently, only very small
    exponents (in terms of spectral radius) need  be exponentiated via
    Lanczos iterations. This enables a far cheaper computation of
    exponentials.
\item 
    Lastly, we introduce a technique to lower the number of exponentials by
    carefully exploiting the time-symmetry of Magnus expansions. This
    involves starting from a slightly modified version of the Magnus
     expansions of \cite{IKS18sinum}, which will be expressed solely in odd
    powers of the time step. %
%
\end{enumerate}
%
%

\subsection{Magnus--Zassenhaus splittings of order two and four}
The simplest method in the family of methods introduced in this paper is the
splitting
\begin{equation}
 \label{eq:expmidpoint}
 \MM{u}^1 = \exp\!\left(\Frac12 \ii \dt \ve \K_2\right) \exp\!\left(- \ii \dt \ve^{-1} \DDD_{\momO} \right) \exp\!\left(\Frac12 \ii \dt \ve \K_2\right) \MM{u}^0.
\end{equation}
where, under spectral collocation as a choice of spatial discretisation,
$\K_k$ is a discretisation of $\dx^k$, $\DDD_f$ is the diagonal matrix with
values of $f$ along the diagonal, and $\momO = \int_0^\dt V(\xi) d\xi$.
The approximation of $\momO$ at the midpoint by $\dt V(\dt/2)$, reminiscent of the familiar {\em Strang splitting,\/} results in the
so-called  {\em exponential midpoint rule\/} with $\O{h^3}$ local
error\footnote{In the context of our Magnus expansions, which are {\em
power-truncated}, time-symmetric and consequently feature only odd powers in
the Taylor expansion, any $\O{h^{2k}}$ quadrature method automatically
becomes $\O{h^{2k+1}}$. In particular, evaluation at middle of the interval
incurs $\O{h^3}$ error. See \cite{IKS18sinum} for a more detailed
narrative.}, which is
well known and has been in use for a long while \cite{lubich08fqc}.\\
%

The first non-trivial splitting in the family of Magnus--Zassenhaus
splittings introduced in this paper is the order four splitting which will be
derived in Section~\ref{sec:MZh4},
\begin{equation}
\label{eq:Zass2}
	 \MM{u}^1 = \ee^{\frac12 \ii  \dt \ve \K_2} \ee^{-\frac12 \ii  \ve^{-1}  \DDD_{\momO}} \ee^{\WW}  \ee^{-\frac12 \ii  \ve^{-1}  \DDD_{\momO}} \ee^{\frac12 \ii  \dt \ve \K_2} \MM{u}^0,
\end{equation}
where $\mom{1}{1} = \int_0^\dt (\xi - h/2) V(\xi) d\xi$ and the central
exponent is
\begin{Eqnarray*}
\WW & = & \Frac{1}{6}\ii  \dt \ve^{-1} \DDD_{\left(\dx \momO \right)^2}
+\left(\K_1 \DDD_{\dx \mom{1}{1}}+ \DDD_{\dx
\mom{1}{1}}\K_1\right)\\
&&-\Frac{1}{3}\ii \dt^2 \ve \left(\K_2 \DDD_{\dx^{2} \momO } + \DDD_{\dx^{2}
\momO } \K_2\right).
\end{Eqnarray*}

\subsection{Organisation of the paper}
Section~2 is devoted to the derivation of our fourth-order scheme
\R{eq:Zass2}. The derivation commences in Subsection~2.1 with a
straightforward rewriting of the fourth-order integral-preserving
simplified-commutator Magnus expansion of \cite{IKS18sinum} in the
semiclassical regime. An analysis of the size of various terms is carried out
in Subsection~2.2, while keeping discretisation considerations in mind,
followed by a symmetric Zassenhaus splitting in Subsection~2.3.
Implementation details for this scheme are discussed briefly in
Subsection~2.4.

The derivation of the sixth-order Magnus--Zassenhaus scheme of Section~3 is
more involved than the approach in Section~2. In contrast to the fourth-order
Magnus expansion, the sixth and higher order Magnus expansions of
\cite{IKS18sinum} feature even powers of the time step. In Subsection~3.2 we
discuss how these even powers result in a large number of exponential stages.
An approach for remedying this issue, which involves expressing the Magnus
expansion solely in odd powers of the time step, is developed in
Subsections~3.3--3.6. We conclude the derivation by performing a Zassenhaus
splitting on such a modified Magnus expansion in Subsection~3.7. 

Numerical examples confirming the efficiency of our schemes are provided in
Section 4, while in the last section we briefly summarise the main
conclusions.

The rules for simplifying commutators from \cite{psalgebra} and the symmetric
Zassenhaus algorithm of \cite{bader14eaf}, which are both utilised in the
derivation of the Zassenhaus splitting in Subsections~2.3 and 3.7, are confined
to Appendices A and B, respectively.

\section{Fourth-order Magnus--Zassenhaus splitting}
\label{sec:MZh4}
Following the approach of \cite{BIKS,IKS18sinum}, the \schr equation
\R{eq:schr} can be written in the form
\begin{equation}\label{eq:Amag}
 u'(t)= A(t) u(t),\qquad t\geq0,
\end{equation}
with $A(t) = \ii [\ve \dx^2 -  \ve^{-1} V(t)]$. Its solution can be
approximated via the Magnus expansion \cite{magnus54}, which is an infinite
series of nested integrals and commutators. In practice we use finite
truncations of this series for a numerical scheme,
\begin{equation}\label{eq:Amag2}
u^{n+1} = \ee^{\Theta_m(t_{n+1},t_n)} u^n,
\end{equation}
where $t_n = t_0 + n\dt$, $\dt$ is the time step, $u^n$ is an approximation
to the exact solution at $t_n$, and $\Theta_m(t,s)$ is a finite truncation of
the Magnus expansion $\Theta(t,s)$ whose exponential evolves the exact solution from $s$
to $t$.

\subsection{Simplified-commutator Magnus expansions}
In the semiclassical regime, the order two and order four
simplified-commutator Magnus expansions of \cite{IKS18sinum} that preserve
the integral become
\begin{Eqnarray*}
\label{eq:Omt1} \Theta_0(h) &=& \ii \dt \ve \dx^2 - \ii  \ve^{-1} \momO,\\
\label{eq:Omt3} \Theta_2(h) &=& \ii \dt \ve \dx^2 - \ii  \ve^{-1} \momO -
2\Ang{1}{\dx \mom{1}{1}},
\end{Eqnarray*}
respectively.  This follows directly from the results of \cite{IKS18sinum}
by accounting for the additional scalar factors  $\ve$ and $\ve^{-1}$. The quantities $\momO$ and $\mom{1}{1}$ have been already  introduced
earlier and
\begin{equation}\label{eq:angbra}
  \ang{f}{k} := \Frac12 \left( f \circ \dx^k + \dx^k \circ f \right), \quad  k \geq 0,\ f \in
  \CC{C}_p^{\infty}(I;\BB{R})
\end{equation}
are the symmetrised differential operators which first appeared in
\cite{bader14eaf}, and have been studied in detail in
\cite{psalgebra,Singh17}. These operators are used here for the same reasons as in \cite{bader14eaf,BIKS,IKS18sinum} -- namely, for
preservation of stability under discretisation after simplification of
commutators.

\subsection{Analysis of size}
As a consequence of \R{eq:dxve}, the symmetrised differential operators scale
as
\begin{equation}\label{eq:angbradisc}
 \ang{f}{k} \leadsto \Frac12 \left( \DDD_f \K_k + \K_k \DDD_f \right) = \O{\ve^{-k}}
\end{equation}
upon discretisation, assuming that $f$ is independent of $\ve$. For
convenience, we abuse notation\footnote{This notion can be made more rigorous
by noting that in the semiclassical regime the derivatives of the solution
grow as $\norm{2}{\dx^k u}\leq C_k \ve^{-k}$, and consequently
$\norm{2}{\ang{f}{k}(u)} \leq C_k \norm{}{f} \ve^{-k}$. The analysis of the
methods presented here can be made rigorous in accordance with these
observations. See \cite{psalgebra,Singh17} for details. } somewhat and use
the shorthand $\ang{f}{k} = \O{\ve^{-k}}$ instead of \R{eq:angbradisc}.

%


\begin{theorem}\label{thm:ve}
All nested commutators of $\ve \dx^2$, $\ve^{-1} V$ and, in general, terms of
the form $\ve^{k} \ang{f}{k+1}$ are $\O{\ve^{-1}}$ or smaller
\cite{psalgebra,Singh17}.
\end{theorem}
Thus a grade $p$ (i.e. $p-1$ times nested) commutator  of $\ii h \ve \dx^2$
and $-\ii h \ve^{-1} V$, which is typically considered to be $\O{h^p}$,
without accounting for $\ve$ or growth corresponding to finer spatial
discretisation, becomes $\O{h^p \ve^{-1}}$ when analysed taking all factors
into account.

To unify the analysis of magnitude of terms, we express the time step in
terms of the inherent parameter $\ve$,
\begin{equation}\label{eq:hve} h = \O{\ve^{\sigma}},
\quad \sigma>0,
\end{equation}
where a smaller $\sigma$ corresponds to larger time steps. We will be aiming
to develop methods that work for small values of $\sigma$ and at the very
least we may assume that $\sigma \leq 1$.

Taylor expansion for an analytic potential, $V(\xi) = \sum_{k=0}^\infty
[V^{(k)}(0)/k!] \xi^k$, shows that $\mom{1}{1}=\int_0^\dt (\xi - h/2) V(\xi)
d\xi$ is $\O{h^3}$ for smooth potentials, since its $\O{h^2}$ term,
$V^{(0)}(0)\int_0^\dt (\xi - h/2) d\xi$, vanishes.

Combining all these ingredients, we conclude that
\begin{equation}
\label{eq:Omt3size}
\Theta_2(h) = \overbrace{\ii \dt \ve \dx^2}^{\O{\ve^{\sigma-1}}} - \overbrace{\ii  \ve^{-1} \momO}^{\O{\ve^{\sigma-1}}} -
\overbrace{2\Ang{1}{\dx \mom{1}{1}}}^{\O{\ve^{3\sigma-1}}}.
\end{equation}
Here it becomes evident that the two leading terms in the Magnus expansion
are significantly larger than the trailing terms. As mentioned in the
introduction, the direct usage of Lanczos iterations for $\exp(\Theta_2)$
 is costly due to the $\O{\ve^{\sigma-1}}$ growth of the
leading terms.

On the other hand, the large term $\ii \dt \ve \dx^2$ can be exponentiated
via FFTs at a $\O{N \log N}$ cost since it discretises as a circulant matrix
and the diagonal term $\ii \ve^{-1} \momO$ can be exponentiated directly at the cost of
$\O{N}$ operations. Thus, we resort to Zassenhaus splittings to separate these
vastly disparate scales, whereby only the small $\O{\ve^{3\sigma -1}}$ terms
need  be exponentiated via Lanczos iterations. The small size of the
remaining exponents means we can work with large time steps (a small
$\sigma$, say $\sigma=1/2$ for instance) and still require very few Lanczos
iterations.

\subsection{Symmetric Zassenhaus splitting of order four}

In order to split the exponential of the fourth-order Magnus expansion,
$\Theta_2$, using Zassenhaus splitting algorithm, we follow the procedure
outlined in Appendix~\ref{sec:Zassenhaus}. In particular, the exponential we
wish to approximate here is $\ee^{\Theta_2}$. Thus, we take
$\mathcal{W}^{[0]}=\Theta_{2}$, which needs to be split in two components, $X
= W^{[0]}$ and $Y = \mathcal{W}^{[0]}- W^{[0]}$. We usually prefer to extract
the largest terms first yet we may choose between $W^{[0]} = \ii \dt \ve
\dx^2$ and $W^{[0]} = -\ii \ve^{-1} \momO$, which are both
$\O{h\ve^{-1}}=\O{\ve^{\sigma-1}}$.


In the specific method presented here, we opt for the former,
whereby $X = W^{[0]}
= \ii \dt \ve \dx^2=\ii \dt \ve \ang{1}{2}$ and $Y = \mathcal{W}^{[0]}-
W^{[0]} =- \ii \ve^{-1} \ang{\momO}{0} - 2 \Ang{1}{\dx \mom{1}{1}}$. (Observe that $\ang{1}{k} = \dx^k$ and $\ang{f}{0} = f$. We
will switch between these  equivalent notations at our convenience.) Using
the sBCH formula \R{eq:sBCH}, we write
\begin{equation}\label{eq:3:sbch_inverted}
  \ee^{\Theta_{2}} = \ee^{\mathcal{W}^{[0]}} =\ee^{\frac12 \ii \dt \ve \dx^2} \ee^{\mathcal{W}^{[1]} } \ee^{\frac12  \ii \ve \dt \dx^2} + \O{h^5 \ve^{-1}},
\end{equation}
where, using the order four expression for $\mathcal{W}^{[1]} =
\CC{sBCH}(-X,X+Y)$ from \R{eq:sBCHshifted},
\begin{Eqnarray}
\label{eq:W1solve}\mathcal{W}^{[1]} &=& - \ii \ve^{-1} \ang{\momO}{0} - 2
\Ang{1}{\dx \mom{1}{1}} +  \Frac{1}{24}  \ii \dt^2 \ve  [[ \ang{\momO}{0} ,\ang{1}{2}], \ang{1}{2}] \\
\nonumber &&+ \Frac{1}{12} \dt^2  \ve [[ \Ang{1}{\dx \mom{1}{1}},\ang{1}{2}],
\ang{1}{2}]
- \Frac{1}{12} \ii \dt   \ve^{-1} [[ \ang{\momO}{0}, \ang{1}{2}],  \ang{\momO}{0}]\\
\nonumber &&- \Frac{1}{6} \dt [[\ang{\momO}{0}, \ang{1}{2}],\Ang{1}{\dx
\mom{1}{1}}] - \Frac{1}{6}\dt  [[
\Ang{1}{\dx \mom{1}{1}}, \ang{1}{2}], \ang{\momO}{0}]\\
\nonumber && + \Frac{1}{3}\ii  \dt \ve   [[ \Ang{1}{\dx \mom{1}{1}},
\ang{1}{2}], \Ang{1}{\dx \mom{1}{1}}] + \O{h^5 \ve^{-1}},
\end{Eqnarray}
and we have ignored all $\O{h^5\ve^{-1}}=\O{\ve^{5\sigma-1}}$ terms in the
sBCH since we are applying the Zassenhaus splitting to the fourth-order
Magnus expansion.


Further, we observe that
\begin{Eqnarray}
\label{eq:W1solve2}\mathcal{W}^{[1]} &=& - \ii \ve^{-1} \ang{\momO}{0} - 2
\Ang{1}{\dx \mom{1}{1}} +  \Frac{1}{24}  \ii \dt^2 \ve  [[ \ang{\momO}{0} ,\ang{1}{2}], \ang{1}{2}] \\
\nonumber && - \Frac{1}{12} \ii \dt   \ve^{-1} [[ \ang{\momO}{0},
\ang{1}{2}],  \ang{\momO}{0}] + \O{h^5 \ve^{-1}},
\end{Eqnarray}
since most commutators are $\O{h^5 \ve^{-1}}$ or smaller due to
Theorem~\ref{thm:ve} and the observation that $\momO = \O{h}$ and $\mom{1}{1}
= \O{h^3}$.
These nested commutators can be simplified using the identities \R{eq:id},
the detailed theory for which is presented in \cite{psalgebra} and which were
utilised in the derivation of the simplified-commutator Magnus expansion
$\Theta_2$ in \cite{IKS18sinum}. Applying these identities, we find
\begin{Eqnarray}
\label{eq:commutators} [\ang{\momO}{0} ,\ang{1}{2}] & = & - 2 \ang{\dx \momO}{1}, \\
\nonumber\left[ \left[\ang{\momO}{0} ,\ang{1}{2}\right], \ang{1}{2}\right] & = &  4 \ang{\dx^2 \momO}{2} -  \ang{\dx^4 \momO}{0},\\
\nonumber \left[ \Ang{1}{\dx \mom{1}{1}},\ang{1}{2}\right] &= & - 2\ang{\dx^2  \mom{1}{1} }{2} + \Frac12 \ang{\dx^4 \mom{1}{1}}{0},\\
\nonumber \left[ \left[ \ang{\momO}{0}, \ang{1}{2}\right],
\ang{\momO}{0}\right] & = &  -2 \ang{(\dx \momO)^2}{0}.
\end{Eqnarray}
Thus, the central exponent in \R{eq:3:sbch_inverted} up to $\O{h^5
\ve^{-1}}=\O{\ve^{5\sigma-1}}$ is
\begin{Eqnarray}\label{eq:3:central_exponent}
\mathcal{W}^{[1]} &=&  - \ii \ve^{-1}  \ang{\momO}{0} +  \Frac{1}{6} \ii
\dt^2 \ve \ang{\dx^2
\momO}{2} - 2 \Ang{1}{\dx \mom{1}{1}}\\
\nonumber && -  \Frac{1}{24} \ii \dt^2  \ve \ang{\dx^4 \momO}{0}+\Frac{1}{6}
\ii \dt \ve^{-1}
 \ang{(\dx \momO)^2}{0}.
\end{Eqnarray}
In order to approximate the exponential of $\mathcal{W}^{[1]}$ (which is
required in \R{eq:3:sbch_inverted}), we apply the sBCH again along the lines
of \R{eq:3:sbch_inverted}, writing
\begin{equation}\label{eq:4:sbch_inverted}
  \ee^{\mathcal{W}^{[1]}} =\ee^{-\frac12  \ii \ve^{-1} \momO} \ee^{\mathcal{W}^{[2]} } \ee^{-\frac12  \ii \ve^{-1} \momO} + \O{h^5\ve^{-1}}.
\end{equation}
Here we have chosen to extract the largest term $X = W^{[1]} = - \ii \ve^{-1}
\momO = \O{h\ve^{-1}}$ and the remainder is $Y = \mathcal{W}^{[1]}- W^{[1]} =
\O{h^3\ve^{-1}}$. The choice of X is more obvious in this case since there is
only one term of size $\O{h\ve^{-1}}$. Keeping the sizes of $X$ and $Y$ in
mind, we find that
\begin{Eqnarray} \label{eq:4:central_exponent}
\nonumber \mathcal{W}^{[2]} &:=& \CC{sBCH}(-X,X+Y)+ \O{h^5 \ve^{-1}} = (-X + (X+Y)) = Y + \O{h^5 \ve^{-1}}\\
&=& \Frac{1}{6} \ii \dt^2 \ve \ang{\dx^2 \momO}{2} - 2 \Ang{1}{\dx
\mom{1}{1}}\\
\nonumber && \mbox{}-  \Frac{1}{24}\ii \dt^2 \ve \ang{\dx^4 \momO}{0}+\Frac{1}{6}
\ii \dt \ve^{-1}
 \ang{(\dx \momO)^2}{0} + \O{h^5 \ve^{-1}}
\end{Eqnarray}
suffices up to $\O{h^5 \ve^{-1}}$ since $[[Y,X],X]$ and $[[Y,X],Y]$ are
$\O{h^5\ve^{-1}}$ and $\O{h^7\ve^{-1}}$, respectively. Further, we discard
the $\O{h^3 \ve}$ term $-  \Frac{1}{24}\ii \dt^2 \ve \ang{\dx^4 \momO}{0}$
from \R{eq:4:central_exponent} since, under the assumption $\sigma \leq 1$,
it is effectively $\O{\ve^{5\sigma-1}}$ or smaller.

Finally, we complete the derivation of our fourth-order Magnus--Zassenhaus
splitting by combining \R{eq:3:sbch_inverted} and \R{eq:4:sbch_inverted},
\begin{equation}
\label{eq:MZ2}
	 \exp(\Theta_{2}) = \ee^{\frac12 W^{[0]}} \ee^{\frac12 W^{[1]}}
									 \ee^{\WW^{[2]}} \ee^{\frac12 W^{[1]}} \ee^{\frac12 W^{[0]}} + \O{\ve^{5\sigma-1}},
\end{equation}
which has the same form as the splittings described in
Appendix~\ref{sec:Zassenhaus}. Here the exponents (chosen and/or computed
during the derivation) are
\begin{Eqnarray}
W^{[0]} & = & \ii \dt \ve \dx^2  = \O{\ve^{\sigma-1}},\label{eq:MZ2W} \\
W^{[1]} & = & - \ii  \ve^{-1}  \momO = \O{\ve^{\sigma-1}}, \nonumber\\
\WW^{[2]} & = & \Frac{1}{6}\ii  \dt \ve^{-1} \left(\dx \momO \right)^2
-2\Ang{1}{\dx \mom{1}{1}}+\Frac{1}{6}\ii  \dt^2 \ve \Ang{2}{\dx^{2} \momO } =
\O{\ve^{3\sigma-1}}.\nonumber
\end{Eqnarray}

It is only by this stage -- having arrived at an asymptotic splitting
expressed in operatorial terms -- that we start considering discretisation
issues. Apart from analysis of size, considerations of spatial and temporal
discretisation (in the form of approximation through quadrature) are entirely
independent of the procedure leading to \R{eq:MZ2W} and of each other: one may proceed to address them in any order.

\begin{remark}
Note that, in contrast to the Magnus--Zassenhaus splitting of \cite{BIKS},
our exponents retain the integrals intact (in $\momO$ and $\mom{1}{1}$). This
provides us with the flexibility to choose the most suitable approximation method
 including the possibility of analytic approximations. This
flexibility proves particularly helpful in the case of potentials with high
temporal oscillations.
\end{remark}

\begin{remark}
As a quick sanity check, we do arrive at the standard symmetric Zassenhaus
splitting \cite{bader14eaf} for time-independent potentials, $V(x,t) :=
V(x)$.
\end{remark}

\subsection{Implementation}

\label{sec:MZ2implementation} The operatorial splitting given by \R{eq:MZ2}
and \R{eq:MZ2W} can be converted into a concrete numerical scheme by
resorting to spatial discretisation. As typical, we use spectral collocation,
whereby $W^{[0]}$ is discretised as the circulant matrix $\widetilde{W^{[0]}}
= \ii \dt \ve \K_2$. Note that $\K_2$ is diagonalisable via Fourier transform
since
\begin{equation}
\label{eq:KkFFT}
\K_k = \FFF^{-1} \DDD_{c_k} \FFF,
\end{equation}
where $c_k$ is the
symbol of $\K_k$, $\FFF$ is the Fast Fourier Transform (FFT) and $\FFF^{-1}$
is its inverse (IFFT). Thus, we can evaluate
\[
\ee^{\frac12 \widetilde{W^{[0]}}} \MM{u}  = \mathcal{F}^{-1} \DDD_{\exp(\frac12 \ii \dt \ve c_2)} \mathcal{F} \MM{u},
\]
in $\O{N \log N}$ operations. The next exponent, $W^{[1]}$, is discretised as
the diagonal matrix $\widetilde{W^{[1]}} = - \ii \ve^{-1} \DDD_{\momO}$ and
exponentiated pointwise in $\O{N}$ operations,
\[
\ee^{\frac12 \widetilde{W^{[1]}}} \MM{u} = \DDD_{\exp(\frac12 \ii \ve^{-1} \momO)} \MM{u}.
\]
Lastly, the innermost exponent,
\begin{Eqnarray*}
\WW^{[2]} \leadsto  \widetilde{\WW^{[2]}} & = & \Frac{1}{6}\ii  \dt \ve^{-1}
\DDD_{\left(\dx \momO \right)^2} +\left(\K_1 \DDD_{\dx \mom{1}{1}}+ \DDD_{\dx
\mom{1}{1}}\K_1\right)\\
&&-\Frac{1}{3}\ii \dt^2 \ve \left(\K_2 \DDD_{\dx^{2} \momO } + \DDD_{\dx^{2}
\momO } \K_2\right),
\end{Eqnarray*}
which is neither circulant nor diagonal, can be very effectively
exponentiated via Lanczos iterations because (unlike the first two exponents)
it is very small in size, even when combined with large time steps.

For instance, when using time steps as large as $h = \O{\ve^{1/2}}$, i.e.\
under the scaling $\sigma = 1/2$, the exponent $\widetilde{\WW^{[2]}}$ is
$\O{\ve^{1/2}}$. As the semiclassical parameter $\ve \ll 1$ becomes small,
this spectral radius approaches $0$. Thus the superlinear accuracy of Lanczos
iterations is seen immediately and merely four Lanczos iterations are required
for an exponentiation to an accuracy of $\O{\ve^{5\sigma-1}}=\O{\ve^{3/2}}$,
which is the accuracy of the splitting. A more detailed analysis of the
number of Lanczos iterations required for such exponents is carried out in
\cite{bader14eaf,BIKS}.
%

Each Lanczos iteration requires the computation of a matrix--vector product
of the form $\widetilde{\WW^{[2]}} \MM{v}$. In general Lanczos iterations for
the discretisation of $W = \sum_{k=0}^n \ii^{k+1} \ang{f_k}{k}$ to
$\widetilde{W} = \ii \DDD_{f_0} + \frac12 \sum_{k=1}^n \ii^{k+1}
\left(\DDD_{f_k} \K_k +\K_k \DDD_{f_k} \right)$ can be efficiently computed
in $2n+2$ Fast Fourier Transforms (included inverses) using the expression
 \begin{Eqnarray}
\label{eq:Wv}  \widetilde{W}\MM{v} = \ii \DDD_{f_0} \MM{v} + \Frac12 \left(
\sum_{k=1}^n \ii^{k+1} \DDD_{f_k} \FFF^{-1} \DDD_{c_k} \right) \FFF \MM{v} +
\Frac12 \FFF^{-1} \left( \sum_{k=1}^n \ii^{k+1}  \DDD_{c_k} \FFF \DDD_{f_k}
\MM{v}\right).
\end{Eqnarray}
Thus, the computation of $\widetilde{\WW^{[2]}} \MM{v}$ in each Lanczos
iteration requires six FFTs and the total cost of its exponentiation is 24
FFTs. This completes the derivation of the splitting described in
\R{eq:Zass2}.

\begin{remark}
$\Theta_0$ is an order two Magnus expansion. An order two
Zassenhaus splitting is merely the Strang splitting. Thus an order two
Magnus--Zassenhaus splitting is, trivially, \R{eq:expmidpoint}, which has been
presented in the introduction.
\end{remark}

\begin{remark}
We assume that the integrals $\momO$ and $\mom{1}{1}$ are either available
analytically or can be adequately approximated via quadrature formulae such
as Gauss--Legendre quadratures, along the lines of \cite{IKS18sinum}.

\end{remark}

\begin{remark} Along the lines of \cite{BIKS,IKS18sinum}, we also assume that spatial derivatives of the potential, as well as
its integrals, are either available or inexpensive to compute.
\end{remark}

\section{Sixth-order Magnus--Zassenhaus splitting}
\label{sec:CFModd}
The order six simplified-commutator Magnus expansion of \cite{IKS18sinum} in
the context of the semiclassical regime, is
\begin{Eqnarray}\label{eq:Th4}
\nonumber \Theta_4(h) & = &  \overbrace{\ii \dt \ve \dx^2 - \ii \ve^{-1} \momO}^{\O{\dt \ve^{-1}}} -  \overbrace{2 \Ang{1}{\dx \mom{1}{1}}}^{\O{\dt^3 \ve^{-1}}} + \overbrace{\ii \ve^{-1} \Lf{\psi}{1}{1} + 2\ii \ve \Ang{2}{\dx^2 \mom{2}{1}}}^{\O{\dt^4 \ve^{-1}}} \\
\nonumber &&\mbox{} +\overbrace{\Frac{1}{6} \Ang{1}{\Lf{\varphi_1}{1}{2} + \Lf{\varphi_2}{2}{1} }}^{\O{\dt^4 \ve^{-1}}} + \overbrace{\Frac{1}{6} \Ang{1}{\Lf{\phi_1}{1}{2} + \Lf{\phi_2}{2}{1} }}^{\O{\dt^5 \ve^{-1}}}\\
\label{eq:Omt5}&&\mbox{}+\overbrace{\Frac{4}{3} \ve^{2} \Ang{3}{\dx^3
\mom{3}{1}}}^{\O{\dt^5 \ve^{-1}}} +\overbrace{\Frac{1}{4} \ii \ve \dx^4
\mom{2}{1}}^{\O{\dt^4 \ve}} = \mTh{} + \O{\dt^7 \ve^{-1}},
\end{Eqnarray}
where $\mom{j}{k}$ are integrals on the line,
\begin{equation}
\label{eq:mom}
\mom{j}{k} = \Int{\zeta}{0}{\dt}{\tilde{B}_{j}^k(\dt,\zeta) V(\zeta) },
\end{equation}
$\Lf{f}{a}{b}$ are integrals over a triangle,
\begin{equation}
\label{eq:Lf1}
\Lf{f}{a}{b} = \Int{\zeta}{0}{\dt}{\Int{\xi}{0}{\zeta}{f(\dt,\zeta,\xi) \left[ \dx^a V(\zeta)\right] \left[\dx^b V(\xi) \right] }},
\end{equation}
$\tilde{B}_j(\dt,\zeta) = \dt^j B_j\left(\zeta/\dt \right)$ are rescaled
Bernoulli polynomials \cite{abram64spfunc,lehmer88ber},
and
\begin{Eqnarray}
\nonumber\psi(\dt,\zeta,\xi) &=& \zeta - \xi - \Frac{\dt}{3},\\
\nonumber\varphi_1(\dt,\zeta,\xi)& = & \dt^2 - 4 \dt \xi + 2 \zeta \xi, \\
\nonumber \varphi_2(\dt,\zeta,\xi)& = & (\dt-2\zeta)^2 - 2 \zeta \xi,\\
\nonumber \phi_1(\dt,\zeta,\xi) & = & \dt^2 - 6 \dt \zeta + 6 \dt \xi + 6 \zeta \xi + 3 \zeta^2 - 12 \xi^2,\\
\nonumber \phi_2(\dt,\zeta,\xi) & = & \dt^2 - 6 \dt \zeta + 6 \dt \xi - 6
\zeta \xi + 5 \zeta^2.
\end{Eqnarray}

\begin{remark}\label{rmk:muindices} Note that, since $B_0(x) = 1$,
\[ \mom{0}{k} = \mom{j}{0} = \mom{0}{0},\quad  j,k\in \BB{Z}^+.\]
To avoid confusion, therefore, such terms will always be indexed as
$\mom{0}{0}$.
\end{remark}

\subsection{Size analysis}

Size of the terms is analysed with the aid of Lemmas~\ref{lem:powersofh_mu}
and \ref{lem:powersofh_lambda}.

\begin{lemma} \label{lem:powersofh_mu}
    For an analytic potential $V$,
    \begin{equation}\label{eq:powersofh_mu}
    \mom{j}{k} = \O{\dt^{jk+1}}, \quad k\neq 1, \qquad \mom{j}{1} = \O{h^{j+2}}.
    \end{equation}

    \begin{proof} \label{pf:powersofh_mu}
    Since the $j${\em th} rescaled Bernoulli polynomial scales as $\O{\dt^j}$, we
    expect $\mom{j}{k} = \O{\dt^{jk+1}}$ . The term $\mom{j}{1}$ gains an extra power of $\dt$ and is
    $\O{\dt^{j+2}}$ since  integrals of  Bernoulli   polynomials vanish,
    $\Int{\zeta}{0}{\dt}{\tilde{B}_j(\dt,\zeta) }= 0$. To see this, consider the Taylor expansion,
    $ V(\zeta) = \sum_{n=0}^\infty [V^{(n)}/n!] \zeta^n$,
    and note that in the expansion for $\mom{j}{1}$,
    the leading order term $ \Int{\zeta}{0}{\dt}{V^{(0)}\tilde{B}_j(\dt,\zeta) }= 0$ vanishes.

    \end{proof}
\end{lemma}

\begin{lemma} \label{lem:powersofh_lambda}
    In general, for a polynomial $r_n(\dt,\zeta,\xi)$ featuring only degree-$n$
    terms in $\dt, \zeta$ and $\xi$, the linear (integral) functional is
    \begin{equation} \label{eq:powersofh_lambda1}
    \Lf{r_n}{a}{b}=\O{\dt^{n+2}}.
    \end{equation}
    For a polynomial $p_n$ whose integral over
    the triangle vanishes,
    $\Int{\zeta}{0}{\dt}{\Int{\xi}{0}{\zeta}{p_n(\dt,\zeta,\xi)}}=0$,
    we gain an extra power of $\dt$ and
    \begin{equation} \label{eq:powersofh_lambda1}
    \Lf{p_n}{a}{b} = \O{\dt^{n+3}}.
    \end{equation}
    In particular, this happens for $\psi, \phi_1$ and $\phi_2$, but not for
    $\varphi_1$ and $\varphi_2$.
\end{lemma}

\subsection{Reduction of stages in Zassenhaus splitting}

The expansion \R{eq:Th4} features both odd and even powers of $h$ -- in
particular, it features terms of size $\O{h^4 \ve^{-1}}$. When we commence
the Zassenhaus procedure  directly from $\Theta_4$, the resulting exponential splitting
\begin{equation}
\label{eq:Zass4suboptimal}
 \ee^{\frac12 W^{[0]}} \ee^{\frac12 W^{[1]}} \ee^{\frac12 W^{[2]}} \ee^{\frac12 W^{[3]}}  \ee^{\frac12 W^{[4]}}\ee^{\WW^{[5]}}\ee^{\frac12 W^{[4]}} \ee^{\frac12 W^{[3]}}  \ee^{\frac12 W^{[2]}}\ee^{\frac12 W^{[1]}}\ee^{\frac12 W^{[0]}},
\end{equation}
features exponents $W^{[3]}=\O{h^4 \ve^{-1}}$ and $\WW^{[5]}=\O{h^6
\ve^{-1}}$ and the power of $h$ is even. Similarly, symmetric Zassenhaus
splitting of a higher order Magnus expansion such as $\Theta_8$ also
features $\O{h^4 \ve^{-1}}$, $\O{h^6 \ve^{-1}}$, $\O{h^8 \ve^{-1}}$ and
$\O{h^{10} \ve^{-1}}$ exponents.

In contrast, the Zassenhaus splittings of \cite{bader14eaf} and \cite{BIKS},
only feature exponents in odd powers of $h$: even powers of
$h$ do not appear in \cite{BIKS} since the integrals in the Magnus expansion
are discretised at the outset by Gaussian quadrature. The appearance of
even powers of $h$ is suboptimal because it results in a larger number of
exponentials here.

To remedy this, we note that the results of \cite{IKS18sinum}, from which we
have commenced our analysis, are based on power truncated Magnus expansions.
These methods are known to be odd in $\dt$ around the midpoint of the
interval ($t_{1/2} = t + h/2$) due to time symmetry of the flow
\cite{iserles00lgm,iserles01tsa}.

It should, therefore, be possible to expand $\Theta_m$ solely in odd powers
of $\dt$ for any choice of the potential $V$. A Zassenhaus splitting
commencing from such an odd-powered expansion will never introduce even
powers of $\dt$ since underlying the procedure is a recursive application of
the symmetric BCH which features only odd-grade commutators. Consequently, in
contrast to the eleven exponentials appearing in \R{eq:Zass4suboptimal}, the
Magnus--Zassenhaus splitting \R{eq:Zass4} that we are about to develop later in this
section features just seven exponentials.

\begin{remark} The extra stages appearing in \R{eq:Zass4suboptimal} can also
be seen as a consequence of lack of symmetry of the scheme.
\end{remark}

\subsection{Time symmetry and shifting the origin}

In this section we derive an expression for $\Theta_4$
expanded solely in odd powers of $\dt$. Recall that, in order to advance from
time $t$ to $t+h$ we use the truncated Magnus expansion,
\begin{displaymath}
u(t+\dt) = \ee^{\Theta_m(t+\dt,t)}u(t),
\end{displaymath}
where $\Theta_m(t+\dt,t)$ can be easily recovered from $\Theta_m(\dt,0)$
(shortened to $\Theta_m(\dt)$) by substituting all occurrences of $V(\zeta)$
by $V(t+\zeta)$. Let us define the midpoint of the interval $[t,t+\dt]$ as
$t_{1/2} = t + \dt/2$. Recalling that the expansion is odd about
$t_{1/2}$, we shift the origin to $t_{1/2}$ by defining
\[W(\zeta):= V(t_{1/2}+\zeta),\]
whereby we need to substitute $V(\zeta)$ with $W(\zeta-\dt/2)$ in the
Magnus expansion. Specifically, for methods described here, this only needs
to be done in the terms $\mom{j}{k}$ and $\Lf{f}{a}{b}$,
\begin{Eqnarray}
\label{eq:momN}
\mom{j}{k} &:=& \Int{\zeta}{0}{\dt}{\tilde{B}_{j}^k(\dt,\zeta)\,W\!\left(\zeta-\frac{\dt}{2}\right)\! },\\
\label{eq:LfN}
\Lf{f}{a}{b} &:=& \Int{\zeta}{0}{\dt}{\Int{\xi}{0}{\zeta}{f(\dt,\zeta,\xi) \left[ \dx^a W\!\left(\zeta-\frac{\dt}{2}\right)\ \dx^b W\!\left(\xi-\frac{\dt}{2}\right) \right] \!}}.
\end{Eqnarray}
As we demonstrate in the rest of this section, this substitution allows us to identify and discard even-powered terms in the Magnus expansion.

\begin{remark}
Since we have shifted the origin to $t_{1/2}$, all odd and even components are to be understood with respect to $0$ from this point onwards. This makes identification of the odd components of the Magnus expansion simpler, assuming that the odd and even components of $W$ can be found.
\end{remark}

\subsection{Odd and even components of $\mu$ and $\Lambda$}
\label{sec:oddeven_mu_lambda} Extending the new definitions of $\mu$ and
$\Lambda$ (\ref{eq:momN}, \ref{eq:LfN}), we define
\begin{equation}
\label{eq:momS}
\momS{\star}{j}{k} := \Int{\zeta}{0}{\dt}{\tilde{B}_{j}^k(\dt,\zeta)\,W^\star\!\left(\zeta-\frac{\dt}{2}\right) },\quad \star \in \{e,o\},
\end{equation}
and
\begin{equation}
\label{eq:LfS}
\LfS{f}{\star}{a}{b} := \Int{\zeta}{0}{\dt}{\Int{\xi}{0}{\zeta}{f(\dt,\zeta,\xi) \left[ \dx^a W\left(\zeta-\frac{\dt}{2}\right)\ \dx^b W\left(\xi-\frac{\dt}{2}\right) \right]^\star}},\quad \star \in \{e,o\}.
\end{equation}
The even and odd components of $\mom{j}{k}$ and $\Lf{f}{a}{b}$ can easily be
expressed in terms of these new definitions,
\[
\begin{tabular}{l l}
  $\left(\mom{j}{k}\right)^{o} = \begin{cases} \momS{e}{j}{k} &\mbox{if } jk \mbox{ is even},\\
\momS{o}{j}{k} &\mbox{if } jk \mbox{ is odd}, \end{cases}$
& \quad$\left(\mom{j}{k}\right)^{e}  =  \begin{cases} \momS{o}{j}{k} &\mbox{if } jk \mbox{ is even},\\
\momS{e}{j}{k} &\mbox{if } jk \mbox{ is odd}, \end{cases}$
\end{tabular} \]
where we exploit the fact that $\tilde{B}_j (\dt, \zeta)$ is odd in $(\dt,
\zeta)$ for odd values of $j$ and even for even values. For instance,
$\momO^o = \momOS{e} =
\Int{\zeta}{0}{\dt}{W^e\!\left(\zeta-\frac{\dt}{2}\right)}$. Note carefully
the difference of usage in subscript and superscript. The odd and even parts
of $\Lf{f}{a}{b}$ are
\begin{Eqnarray*}
\left(\Lf{f}{a}{b}\right)^{\!o} &=&  \LfS{f^o}{e}{a}{b}+ \LfS{f^e}{o}{a}{b},\\
\left(\Lf{f}{a}{b}\right)^{\!e} &=&  \LfS{f^o}{o}{a}{b}+ \LfS{f^e}{e}{a}{b}.
\end{Eqnarray*}
For an odd function such as $\psi(\dt,\zeta,\xi): = \zeta - \xi -
\Frac{\dt}{3}$, for instance, this reduces to
\[\left(\Lf{\psi}{a}{b}\right)^{\!o}  = \LfS{\psi}{e}{a}{b}. \]

\subsection{Expanding $\mu$ and $\Lambda$ in powers of $h$.}

\label{sec:powersofh_oddeven_mu_lambda} A natural consequence of separating
odd and even components is that the sizes of certain $\mu$ and $\Lambda$
integrals become smaller (i.e.\ they are associated with higher powers of $h$) than prescribed by
Lemmas~\ref{lem:powersofh_mu} and \ref{lem:powersofh_lambda}. Essentially, this
is due to certain terms disappearing in the relevant Taylor expansions.
\begin{lemma}\label{lem:powersofh_oddeven_mu}
    For an analytic potential $V$  (and thus $W$),
    \begin{equation} \label{eq:powersofh_even_mu}
        \momS{e}{j}{k} = \O{h^{jk+1}}, \quad k \neq 1, \qquad \momS{e}{j}{1} = \O{h^{j+3}},
    \end{equation}
    and
    \begin{equation} \label{eq:powersofh_odd_mu}
        \momS{o}{j}{k} = \O{h^{jk+2}}, \quad k \geq 0.
    \end{equation}

    \begin{proof}
    \label{pf:powersofh_oddeven_mu} For an analytic potential,
    $W\left(\zeta-\dt/2\right) = \sum_{n=0}^{\infty} [W^{(n)}(0)/n!] \left(\zeta-\dt/2\right)^n $,
    we expect a gain of a single power of $\dt$ in $\mom{j}{1}$ since
    $\Int{\zeta}{0}{\dt}{ \tilde{B}_j(\dt,\zeta) W(0)}$ vanishes (see
    Lemma~\ref{lem:powersofh_mu}). Thus, instead of being $\O{\dt^{j+1}}$, this
    term turns out to be $\O{\dt^{j+2}}$.

    In the case of $\momS{e}{j}{1}$ where we take the even part of the potential
    $W$ at $\zeta-\frac{\dt}{2}$, we see a gain of two powers of $\dt$ since the
    smallest non-vanishing term is
    \[\frac12 \Int{\zeta}{0}{\dt}{ B_j(\dt,\zeta) W^{(2)}(0) \left(\zeta-\frac{\dt}{2}\right)^2}.\]
    Thus, $\momS{e}{j}{1} = \O{h^{j+3}}$ instead of being $\O{h^{j+2}}$. This is
    the only exception to Lemma~\ref{lem:powersofh_mu} for $\momS{e}{j}{k}$.
    Similarly, for $\momS{o}{j}{k}$, observe that there is no $W^{(0)}$ term
    in the Taylor expansion.
    \end{proof}
\end{lemma}

\begin{lemma} \label{lem:powersofh_oddeven_lambda}
For a polynomial $r_n(\dt,\zeta,\xi)$ featuring only degree-$n$ terms in
$\dt, \zeta$ and $\xi$,
\begin{equation} \label{eq:powersofh_even_lambda_r}
    \LfS{r_n}{e}{1}{1} = \O{h^{n+2}},
\end{equation}
while for a polynomial $p_n$ whose integral over the triangle vanishes,
\begin{equation*}
\label{eq:IntPhi}
 \Int{\zeta}{0}{\dt}{\Int{\xi}{0}{\zeta}{p_n(\dt,\zeta,\xi)}}=0,
\end{equation*}
\begin{equation} \label{eq:powersofh_even_lambda_p}
    \LfS{p_n}{e}{1}{1} = \O{h^{n+4}}.
\end{equation}
In particular, the exception \R{eq:powersofh_even_lambda_p} holds for $\psi,
\phi_1$ and $\phi_2$, but not for $\varphi_1$ and $\varphi_2$. For the odd
part,

\begin{equation} \label{eq:powersofh_even_lambda_r}
\LfS{r_n}{o}{1}{1} = \O{h^{n+3}}, \end{equation}
whether the integral of $r_n$ vanishes or not.
\end{lemma}
The proof for Lemma~\ref{lem:powersofh_oddeven_lambda} follows directly
along the lines of Lemma~\ref{lem:powersofh_oddeven_mu}.

\subsection{Identifying the odd components of Magnus expansions}
\label{sec:extracting_odd_Magnus} Recall that the time symmetry of power
truncated Magnus expansions ensures that, having shifted the origin to $t_{1/2}$, $\Theta_m$ is odd about the origin.
Therefore the even part
$\Theta_{m}^{e}$ vanishes, and it suffices to work with $\Theta_{m} =
\Theta_{m}^{o}$.
For the sixth-order simplified-commutator Magnus expansion of
\cite{IKS18sinum}, $\Theta_4$, we can identify the odd components using the
results of Subsection~\ref{sec:oddeven_mu_lambda},
\begin{Eqnarray}
\label{eq:Theta_4o}
\Theta_4 = \Theta_4^o & = &  \overbrace{\ii \dt \ve \dx^2 - \ii \ve^{-1} \momOS{e}}^{\O{h \ve^{-1}}} -  \overbrace{2 \Ang{1}{\dx \momS{o}{1}{1}}}^{\O{h^3 \ve^{-1}}}\\
\nonumber &&\mbox{} + \overbrace{\ii \ve^{-1} \LfS{\psi}{e}{1}{1} + 2\ii \ve \Ang{2}{\dx^2 \momS{e}{2}{1}}}^{\O{h^5 \ve^{-1}}} \\
\nonumber &&\mbox{} +\overbrace{\Frac{1}{6}\ve^{-1} \Ang{1}{\LfS{\varphi_1 + \phi_1}{o}{1}{2} + \LfS{\varphi_2 + \phi_2}{o}{2}{1} }}^{\O{h^5 \ve^{-1}}} \\
\nonumber &&\mbox{}+\overbrace{\Frac{4}{3} \ve^2 \Ang{3}{\dx^3
\momS{o}{3}{1}}}^{\O{h^5 \ve^{-1}}} +\overbrace{\Frac{1}{4} \ii \ve \dx^4
\momS{e}{2}{1}}^{\O{h^5 \ve}} = \mTh{} + \O{h^7 \ve^{-1}},
\end{Eqnarray}
where size analysis in powers of $h$ is done using
Lemmas~\ref{lem:powersofh_oddeven_mu} and \ref{lem:powersofh_oddeven_lambda},
and we have used the fact that $\left(\Lf{\varphi_j + \phi_j}{a}{b}\right)^o
= \LfS{\varphi_j + \phi_j}{o}{a}{b}$, since $\varphi_j + \phi_j$ is even. We
can now discard the $\O{h^5 \ve}$ term $\Frac{1}{4} \ii \ve \dx^4
\momS{e}{2}{1}$ since it is $\O{\ve^{7\sigma-1}}$ or smaller under the
assumption $\sigma\leq 1$.

By this stage we have been able to eradicate all even powers of $\dt$ and are
now in a position to effectively apply the Zassenhaus splitting to the
sixth-order Magnus expansion $\Theta_4^o$.
%
%

\subsection{Sixth-order symmetric Zassenhaus splitting}
%
A Zassenhaus splitting commencing from the Magnus expansion $\Theta_4^o$
naturally yields an exponential splitting featuring only $\O{\ve^{\sigma
-1}}$, $\O{\ve^{3\sigma -1}}$ and $\O{\ve^{5\sigma -1}}$ terms,
\begin{equation}
\label{eq:Zass4}
 \exp\left(\Theta_4\right) = \ee^{\frac12 W^{[0]}} \ee^{\frac12 W^{[1]}} \ee^{\frac12 W^{[2]}} \ee^{\WW^{[3]}} \ee^{\frac12 W^{[2]}}\ee^{\frac12 W^{[1]}}\ee^{\frac12 W^{[0]}} + \O{\ve^{7\sigma -1}} ,
\end{equation}
where
\begin{Eqnarray*}
W^{[0]}  &=& \ii  \dt \ve \dx^2=\O{\ve^{\sigma-1}},\\
W^{[1]} &=& -\ii  \ve^{-1}  \momOS{e} =\O{\ve^{\sigma-1}},\\
W^{[2]} &=& \Frac16\ii  \dt \ve^{-1}  \left(\dx \momOS{e}\right)^2   - 2 \Ang{1}{\dx \momS{o}{1}{1}} + \Frac{1}{6}\ii  \dt^{2} \ve \Ang{2}{ \dx^2 \momOS{e} }=\O{\ve^{3\sigma-1}},\\
\WW^{[3]} &=& \ii  \ve^{-1}  \LfS{\psi}{e}{1}{1}  - \Frac{1}{24}\ii  \dt^{2} \ve \left(\dx^{4}  \momOS{e} \right)
-\Frac{1}{6}\ii  \ve^{-1}  \left(\dx \momOS{e}\right)^2 \left(\dx^2 \momS{e}{2}{1}\right)\\
&&\mbox{} - \Frac{2}{45} \ii \dt^2 \ve^{-1} \left(\dx \momOS{e} \right)^2 \left(\dx^2 \momOS{e} \right) \\
&& +\Frac16 \Ang{1}{\LfS{\varphi_1+\phi_1}{o}{1}{2} + \LfS{\varphi_2+\phi_2}{o}{2}{1} }\\
&&\mbox{}- \dt \Ang{1}{\left(\dx \momOS{e}\right) \left(\dx^{2} \momS{o}{1}{1} \right) -\Frac13 \left(\dx^2 \momOS{e}\right) \left(\dx \momS{o}{1}{1} \right)} \\
&&\mbox{}+\Frac{1}{30} \ii \dt^3 \ve \Ang{2}{\left(\dx^2 \momOS{e}\right)^2 - 2\left(\dx \momOS{e}\right) \left(\dx^3 \momOS{e}\right)} \\
&&\mbox{}+2 \ii  \ve\Ang{2}{\dx^{2} \momS{e}{2}{1}} + \Frac{4}{3} \ve^{2}\Ang{3}{\dx^3 \momS{o}{3}{1}} \\
&&\mbox{}+\Frac13 \dt^2 \ve^2 \Ang{3}{\dx^3 \momS{o}{1}{1}} - \Frac{1}{120}
\ii \dt^4 \ve^3 \Ang{4}{\dx^4 \momOS{e}} = \O{\ve^{5\sigma-1}},
\end{Eqnarray*}

The term $\Frac{1}{24}\ii  \dt^{2} \ve \left(\dx^{4}  \momOS{e} \right)$ in
$\WW^{[3]}$ is $\O{\ve^{3\sigma+1}}$. We remind the reader that this
splitting is obtained subject to $\sigma\leq 1$. For $1/2 <\sigma \leq 1$,
this term can be combined with the $\O{\ve^{5\sigma-1}}$ terms in
$\WW^{[3]}$, whereas for $\sigma \leq 1/2$ it can be ignored since it is
smaller than $\O{\ve^{7\sigma-1}}$.

\subsection{Remarks concerning implementation}
The approximation of the exponentials in \R{eq:Zass4} can be done along the
same lines as Subsection~\ref{sec:MZ2implementation}. For instance, for
$h=\O{\ve^{1/2}}$, the exponents $W^{[2]}$ and $\WW^{[3]}$ are
$\O{\ve^{1/2}}$ and $\O{\ve^{3/2}}$, respectively. Approximating their
exponentials to an accuracy of $\O{\ve^{7\sigma-1}} = \O{\ve^{5/2}}$ requires
six and two Lanczos iterations, respectively. Each Lanczos iteration requires
the computation of a matrix--vector product via \R{eq:Wv}.

In  each time step of the order-six scheme we need to approximate five line
integrals
($\momOS{e}$, $\momS{o}{1}{1}$, $\momS{e}{1}{2}$, $\momS{e}{2}{1}$, $\momS{o}{3}{1}$) %
and three integrals over the triangle ($\LfS{\phi}{e}{1}{1}$, $\LfS{\psi_1 +
\theta_1}{o}{1}{2}$ and $\LfS{\psi_2 + \theta_2}{o}{2}{1}$).
If analytic expressions are not available, it is possible to approximate
these through quadrature formulae such as Gauss--Legendre quadrature. For the
$\O{\ve^{7 \sigma -1}}$ splitting here, for instance, we require merely three
Gauss--Legendre knots $t_k = h (1 + k \sqrt{3/5})/2 , k=-1,0,1$, with weights
$w_k = \Frac{5}{18}\dt, \Frac49 \dt, \Frac{5}{18} \dt$ \cite{davis84quad},
respectively, for the $\O{\dt^7}$ accuracy which is required of the
quadrature.

Evaluation for the integrals over the line, $\momS{\star}{j}{k}$, follows
directly using these quadrature knots. Integrals over the triangle are also
easily evaluated using the same three knots by substituting $W$ with its
interpolating polynomial. A more detailed narrative on the approximation of
such integrals can be found in \cite{IKS18sinum}.

\section{Numerical examples}
\label{sec:numerics}
    \label{sec:experiments}
\begin{figure}[tbh]
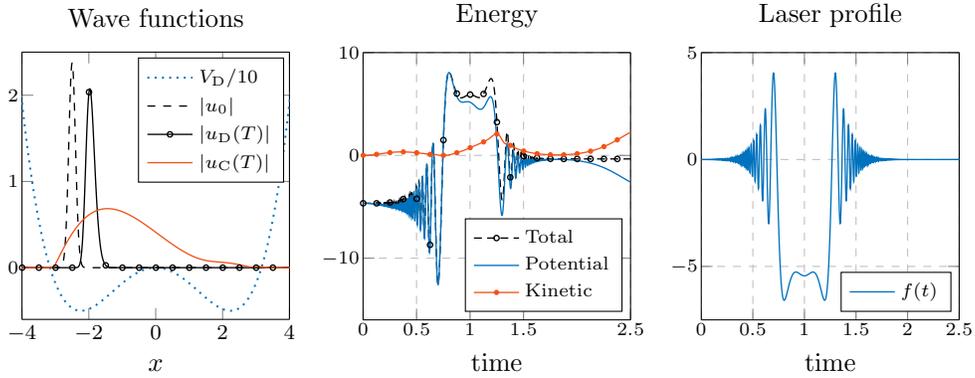

            \input{InitialFinal.tex}
            \input{Energy.tex}
            \input{LaserProfile.tex}
        \caption{The initial condition $u_0$ evolves to $u_{\mathrm{D}}$ under the influence of $V_{\mathrm{D}}$ and
        to $u_{\mathrm{C}}$ under $V_{\mathrm{C}}$ (left); Evolution of energy under $V_{\mathrm{C}}$ (centre); Laser profile $f(t)$ (right).
        Here, the semiclassical parameter is $\ve = 10^{-2}$ and the potential $V_{\mathrm{D}}$ is scaled down for ease of presentation. }
        \label{fig:uIF}
    \end{figure}

The initial condition for our numerical experiments is a Gaussian wavepacket
\[u_0(x) = (\delta  \pi)^{-1/4}  \exp\left( (-(x - x_0)^2 )/ (2 \delta)\right),\quad x_0 = -2.5,\quad \delta = 10^{-2},\]
sitting in the left well of a double well potential,
\[V_{\mathrm{D}}(x) = \Frac15 x^4 - 2 x^2. \]
We take $[-5,5]$ as our spatial domain and $[0,\frac52]$ as our temporal domain.
In the first instance, we consider the behaviour subject to $\ve = 10^{-2}$. When
we allow the wave function $u_0$ to evolve to $u_{\mathrm{D}}(T)$ under the
influence of $V_{\mathrm{D}}$ alone, it remains largely confined to the left
well at the final time, $T=\frac52$ (Figure~\ref{fig:uIF}, left).


{\bf Time-dependent potential.} Superimposing a time dependent excitation of
the form $f(t)\,x$ on the potential -- used for modeling laser interaction --
we are able to exert control on the wave function. The time profile of the
laser used here is
\[ f(t) = 10 \exp(-10 (t - 1)^2) \sin((500 (t - 1)^4 + 10)), \]
which is a highly oscillatory function (Figure~\ref{fig:uIF}, right).
Electric fields of this form, where the instantaneous frequency changes, are
called {\em chirped} pulses. These are used routinely in laser control
\cite{AmstrupChirped}. Even more oscillatory electric fields often result
from optimal control algorithms \cite{MeyerOptimal,CoudertOptimal}.

The effective time-dependent potential is
\[ V_{\mathrm{C}}(x,t) = V_{\mathrm{D}}(x) + f(t) x.\]
Under the influence of the time-dependent potential $V_{\mathrm{C}}(x,t)$,
the initial wave-function $u_0$ evolves to $u_{\mathrm{C}}(T)$, which is not
confined to the left well (Figure~\ref{fig:uIF}, left). In this case, the
total energy is not conserved, stabilising at a higher level once the
external influence vanishes (Figure~\ref{fig:uIF}, centre).

\begin{figure}[tbh]
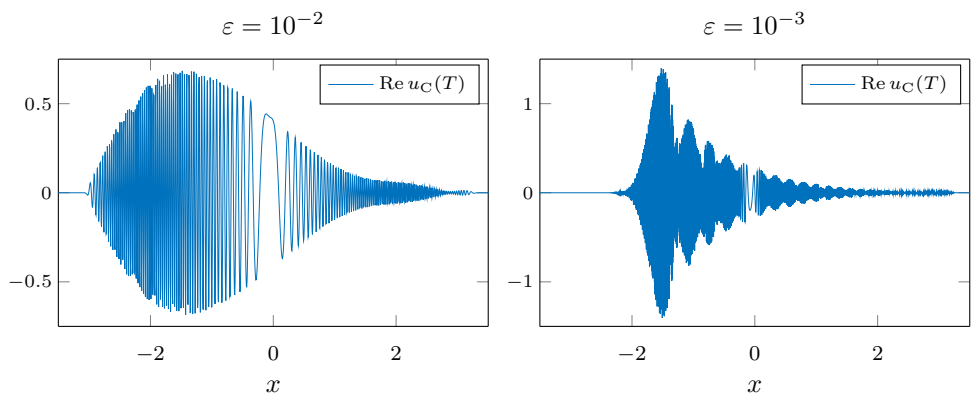

        \input{e1by100oscillations.tex}
        \input{e1by1000oscillations.tex}
    \caption{The highly oscillatory nature of the wavefunction for $\ve=10^{-2}$ and $10^{-3}$ is evident in
    the phase. }
    \label{fig:oscillatory}
\end{figure}


{\bf Methods.} In this section we use IKS to refer to the new
sixth-order method presented in this paper \R{eq:Zass4}. We expect the
asymptotic nature of our splitting \R{eq:Zass4} to outperform methods based
on Lanczos exponentiation of Hamiltonians (or its perturbations) due the
growth of spectral radius in the semiclassical regime. For demonstrating
this, we compare the accuracy and cost of our method with the sixth-order
optimised method CF6:5Opt proposed in \cite{alvermann2011hocm}, labeled as AF
in this section. This is a Lanczos-based method which is highly effective
 in the atomic regime ($\ve=1$). However, it is not specialised for the
semiclassical regime. We also contrast the accuracy of our method to a method
specialised for the semiclassical regime -- the sixth-order method of
\cite{BIKS}, labeled as BIKS in this section. However, this method is not
designed for handling highly oscillatory potentials.

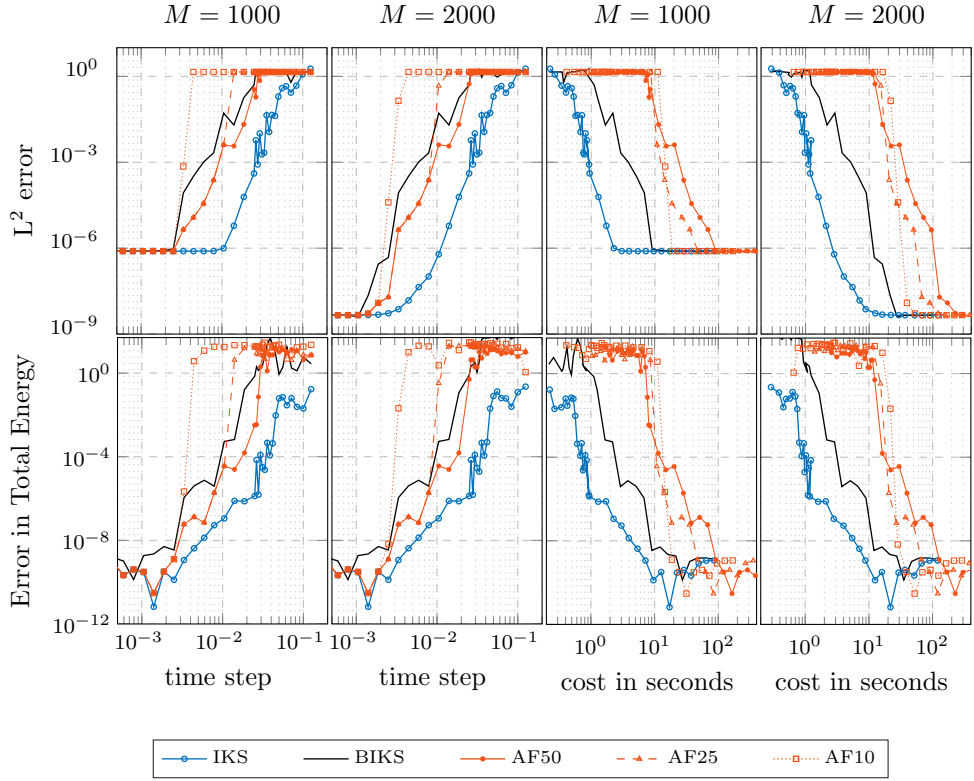
\begin{figure}[tbh]
    \centering
%
%
\definecolor{mycolor1}{rgb}{0.00000,0.44700,0.74100}%
\definecolor{mycolor2}{rgb}{0.85000,0.32500,0.09800}%
\definecolor{mycolor3}{rgb}{0.92900,0.69400,0.12500}%
\definecolor{mycolor4}{rgb}{0.49400,0.18400,0.55600}%
\begin{tikzpicture}

\begin{axis}[%
width=1.1in, height=1.5in, scale only axis, xmode=log, xmin=0.0005, xmax=0.2,
xminorticks=true, ymode=log, ymin=1e-9, ymax=10, yminorticks=true, axis
background/.style={fill=white}, ylabel = {$\CC{L}^2$ error}, yticklabel style
= {xshift=-1ex, yshift=1ex, text height=2ex}, xticklabels={}, grid=both,
title={$M=1000$}]
\addplot [color=colorclassyblue,line width=0.5pt, solid, mark=o,mark
size=1,mark options={solid},mark repeat=1]
  table[row sep=crcr]{%
0.125	1.84227775654152\\
0.1	1.17265319593959\\
0.0833333333333333	0.4732466153733\\
0.0714285714285714	0.27053796669742\\
0.0625	0.450235287063338\\
0.0555555555555556	0.386345109425166\\
0.05	0.19631810855591\\
0.0454545454545455	0.0416870997865452\\
0.0416666666666667	0.04468753649014\\
0.0384615384615385	0.0114665093876507\\
0.0357142857142857	0.0436778345414804\\
0.0333333333333333	0.00214879119160973\\
0.03125	0.00187304850615792\\
0.0294117647058824	0.0100873067378623\\
0.0277777777777778	0.000837317401088285\\
0.0263157894736842	0.00583994127777149\\
0.025	0.000406289736000753\\
0.0186567164179104	6.11838255075089e-05\\
0.0140449438202247	6.03342676943999e-06\\
0.0105042016806723	9.94019024612788e-07\\
0.00788643533123028	7.91513048319724e-07\\
0.00592417061611374	7.87073445620558e-07\\
0.0044404973357016	7.8648980719016e-07\\
0.00333333333333333	7.86692037205271e-07\\
0.0025	7.86858135252198e-07\\
0.00187406296851574	7.86963347129766e-07\\
0.0014052838673412	7.87004767794613e-07\\
0.0010539629005059	7.87039954427015e-07\\
0.000790388871324692	7.87053405105163e-07\\
0.000592838510789661	7.87065032174514e-07\\
0.00044452347083926	7.87070991842819e-07\\
0.000333377783704494	7.87074276806681e-07\\
0.00025	7.87076615494501e-07\\
};

\addplot [color=black, line width=0.5pt, solid]
  table[row sep=crcr]{%
0.125	1.42142337874614\\
0.1	1.41543353882163\\
0.0833333333333333	1.42618001574328\\
0.0714285714285714	0.616961466424723\\
0.0625	1.36502197453387\\
0.0555555555555556	1.31322899363287\\
0.05	1.48065673810313\\
0.0454545454545455	1.40929566131312\\
0.0416666666666667	1.51482248887002\\
0.0384615384615385	1.50566469199867\\
0.0357142857142857	1.61449457806602\\
0.0333333333333333	1.6896145376532\\
0.03125	1.60202745073351\\
0.0294117647058824	1.53472430191142\\
0.0277777777777778	1.6715422990589\\
0.0263157894736842	1.45625161453631\\
0.025	0.500329870104081\\
0.0186567164179104	0.180544404447446\\
0.0140449438202247	0.0199669813715889\\
0.0105042016806723	0.0519286846596421\\
0.00788643533123028	0.0021085765151355\\
0.00592417061611374	0.00101666768566824\\
0.0044404973357016	0.000323378572772979\\
0.00333333333333333	8.54421422486014e-05\\
0.0025	9.18676598854016e-07\\
0.00187406296851574	8.34437282123676e-07\\
0.0014052838673412	7.87312171250678e-07\\
0.0010539629005059	7.8703713763302e-07\\
0.000790388871324692	7.8705293254456e-07\\
0.000592838510789661	7.87065057624439e-07\\
0.00044452347083926	7.87071030196728e-07\\
0.000333377783704494	7.87074294404881e-07\\
0.00025	7.87076620608476e-07\\
};

\addplot [color=colorclassyorange, solid, mark=10-pointed star,mark
size=1,mark options={solid},mark repeat=1]
  table[row sep=crcr]{%
0.125	1.40415419118279\\
0.1	1.40708122746187\\
0.0833333333333333	1.41475609185888\\
0.0714285714285714	1.41466864823918\\
0.0625	1.4143050722728\\
0.0555555555555556	1.41131467275486\\
0.05	1.41438538641492\\
0.0454545454545455	1.41642573597532\\
0.0416666666666667	1.41323061481149\\
0.0384615384615385	1.4142279935137\\
0.0357142857142857	1.4142122496858\\
0.0333333333333333	1.40199193304616\\
0.03125	1.36884212575885\\
0.0294117647058824	0.713325469965504\\
0.0277777777777778	1.13312138961445\\
0.0263157894736842	0.18721245822371\\
0.025	0.351335182829564\\
0.0186567164179104	0.0208632431077521\\
0.0140449438202247	0.00366405373153227\\
0.0105042016806723	0.00405105463605388\\
0.00788643533123028	0.000235577286684763\\
0.00592417061611374	3.60277378172235e-05\\
0.0044404973357016	1.17461176932193e-05\\
0.00333333333333333	4.4586064492763e-06\\
0.0025	7.87305987819122e-07\\
0.00187406296851574	7.87161176196236e-07\\
0.0014052838673412	7.87082965662812e-07\\
0.0010539629005059	7.87078658505635e-07\\
0.000790388871324692	7.87078737212971e-07\\
0.000592838510789661	7.8707858793019e-07\\
0.00044452347083926	7.87078575493616e-07\\
0.000333377783704494	7.87078525006732e-07\\
0.00025	7.87079019076335e-07\\
};

\addplot [color=colorclassyorange, dashed, mark=triangle,mark size=1,mark
options={solid},mark repeat=1]
  table[row sep=crcr]{%
0.125	1.41282517662883\\
0.1	1.41443157385466\\
0.0833333333333333	1.41412757202444\\
0.0714285714285714	1.41096148838448\\
0.0625	1.4144165082393\\
0.0555555555555556	1.41693850510189\\
0.05	1.41437644821367\\
0.0454545454545455	1.41418644260865\\
0.0416666666666667	1.41414128624683\\
0.0384615384615385	1.41295681298473\\
0.0357142857142857	1.41464728915028\\
0.0333333333333333	1.41423784240626\\
0.03125	1.41427256940579\\
0.0294117647058824	1.41421246286528\\
0.0277777777777778	1.41421160656366\\
0.0263157894736842	1.41421223135242\\
0.025	1.41421231095379\\
0.0186567164179104	1.41416788356505\\
0.0140449438202247	1.41421398941181\\
0.0105042016806723	0.0040513327419151\\
0.00788643533123028	0.000235577286690319\\
0.00592417061611374	3.60277378130524e-05\\
0.0044404973357016	1.17461176915773e-05\\
0.00333333333333333	4.45860644940511e-06\\
0.0025	7.87305987497348e-07\\
0.00187406296851574	7.87161176142256e-07\\
0.0014052838673412	7.87082965541765e-07\\
0.0010539629005059	7.87078658567248e-07\\
0.000790388871324692	7.87078736916059e-07\\
0.000592838510789661	7.87078587900957e-07\\
0.00044452347083926	7.87078574762062e-07\\
0.000333377783704494	7.87078526117217e-07\\
0.00025	7.8707901680075e-07\\
};

\addplot [color=colorclassyorange, densely dotted, mark=square,mark
size=1,mark options={solid},mark repeat=1]
  table[row sep=crcr]{%
0.125	1.41422465226304\\
0.1	1.41776764893643\\
0.0833333333333333	1.41380764671022\\
0.0714285714285714	1.41323939994499\\
0.0625	1.41419217613424\\
0.0555555555555556	1.41001414309911\\
0.05	1.40987510823398\\
0.0454545454545455	1.41416218594551\\
0.0416666666666667	1.41562183619228\\
0.0384615384615385	1.41422562242581\\
0.0357142857142857	1.41784190360511\\
0.0333333333333333	1.41421114359345\\
0.03125	1.41421356252877\\
0.0294117647058824	1.41421350441019\\
0.0277777777777778	1.41421347502064\\
0.0263157894736842	1.41421386734806\\
0.025	1.41420604231314\\
0.0186567164179104	1.41421356205934\\
0.0140449438202247	1.41421356237804\\
0.0105042016806723	1.41421356237033\\
0.00788643533123028	1.41421356235043\\
0.00592417061611374	1.4142135623013\\
0.0044404973357016	1.41421566540595\\
0.00333333333333333	0.000709470511975798\\
0.0025	7.87254286335881e-07\\
0.00187406296851574	7.87162463721254e-07\\
0.0014052838673412	7.87082973911847e-07\\
0.0010539629005059	7.87078659209513e-07\\
0.000790388871324692	7.8707873755265e-07\\
0.000592838510789661	7.87078588145297e-07\\
0.00044452347083926	7.87078574197552e-07\\
0.000333377783704494	7.87078525410382e-07\\
0.00025	7.87079017918864e-07\\
};

\end{axis}
\end{tikzpicture}%
        \hspace{-0.38cm}
%
%
\definecolor{mycolor1}{rgb}{0.00000,0.44700,0.74100}%
\definecolor{mycolor2}{rgb}{0.85000,0.32500,0.09800}%
\definecolor{mycolor3}{rgb}{0.92900,0.69400,0.12500}%
\definecolor{mycolor4}{rgb}{0.49400,0.18400,0.55600}%
\begin{tikzpicture}

\begin{axis}[%
width=1.1in, height=1.5in, scale only axis, xmode=log, xmin=0.0005, xmax=0.2,
xminorticks=true, ymode=log, ymin=1e-9, ymax=10, yticklabels={},
yminorticks=true, axis background/.style={fill=white}, xticklabels={},
grid=both, title={$M=2000$}]
\addplot [color=colorclassyblue,line width=0.5pt, solid, mark=o,mark
size=1,mark options={solid},mark repeat=1]
  table[row sep=crcr]{%
0.125	1.85494302749674\\
0.1	1.35892792944903\\
0.0833333333333333	0.488130494494647\\
0.0714285714285714	0.26497196668039\\
0.0625	0.459677280081579\\
0.0555555555555556	0.380603252191244\\
0.05	0.198323808337455\\
0.0454545454545455	0.0524171094127853\\
0.0416666666666667	0.0447008597341876\\
0.0384615384615385	0.0114680381560209\\
0.0357142857142857	0.0436777561208933\\
0.0333333333333333	0.00214763527462076\\
0.03125	0.00187279773846845\\
0.0294117647058824	0.0100874432412646\\
0.0277777777777778	0.000837311912421154\\
0.0263157894736842	0.00583994129926332\\
0.025	0.000406288999606731\\
0.0186567164179104	6.11789073697043e-05\\
0.0140449438202247	5.98275056956967e-06\\
0.0105042016806723	6.11907841941008e-07\\
0.00788643533123028	1.02235240051469e-07\\
0.00592417061611374	4.37777364948667e-08\\
0.0044404973357016	1.49455166558599e-08\\
0.00333333333333333	7.40374031154423e-09\\
0.0025	5.30796338809436e-09\\
0.00187406296851574	4.76959759369982e-09\\
0.0014052838673412	4.63303159501684e-09\\
0.0010539629005059	4.58366429543514e-09\\
0.000790388871324692	4.57348720372152e-09\\
0.000592838510789661	4.54248469051601e-09\\
0.00044452347083926	4.53968159488868e-09\\
0.000333377783704494	4.53222928927586e-09\\
0.00025	4.61340414162414e-09\\
};

\addplot [color=black, line width=0.5pt, solid]
  table[row sep=crcr]{%
0.125	1.4193637600574\\
0.1	1.51409377932297\\
0.0833333333333333	1.44034071060893\\
0.0714285714285714	1.36201421285788\\
0.0625	1.14941806946162\\
0.0555555555555556	0.941130481007499\\
0.05	1.2270336353767\\
0.0454545454545455	1.26519672881492\\
0.0416666666666667	1.33570829005342\\
0.0384615384615385	1.41713762884927\\
0.0357142857142857	0.938311311543404\\
0.0333333333333333	1.77148898716842\\
0.03125	1.6166610234424\\
0.0294117647058824	1.53954904178951\\
0.0277777777777778	1.64889286193774\\
0.0263157894736842	1.47051166807726\\
0.025	0.502688848888131\\
0.0186567164179104	0.1767862286821\\
0.0140449438202247	0.0199669813553296\\
0.0105042016806723	0.0519286846538669\\
0.00788643533123028	0.00210857636874647\\
0.00592417061611374	0.00101666738190946\\
0.0044404973357016	0.000323377616812039\\
0.00333333333333333	8.54385208755247e-05\\
0.0025	4.74182934984492e-07\\
0.00187406296851574	2.77483536405477e-07\\
0.0014052838673412	2.24814075849563e-08\\
0.0010539629005059	4.07126373065946e-09\\
0.000790388871324692	4.49145664427237e-09\\
0.000592838510789661	4.5468056509284e-09\\
0.00044452347083926	4.54614140448793e-09\\
0.000333377783704494	4.53513857704442e-09\\
0.00025	4.61449291227918e-09\\
};

\addplot [color=colorclassyorange, solid, mark=10-pointed star,mark
size=1,mark options={solid},mark repeat=1]
  table[row sep=crcr]{%
0.125	1.4295074309552\\
0.1	1.40729746572539\\
0.0833333333333333	1.41721791013935\\
0.0714285714285714	1.41240231872207\\
0.0625	1.41526693378122\\
0.0555555555555556	1.4142365838197\\
0.05	1.41439238165919\\
0.0454545454545455	1.41419980934031\\
0.0416666666666667	1.41427082254659\\
0.0384615384615385	1.41421400623906\\
0.0357142857142857	1.4142201367728\\
0.0333333333333333	1.41418020835033\\
0.03125	1.41421348023915\\
0.0294117647058824	1.41422876944897\\
0.0277777777777778	1.40992447977729\\
0.0263157894736842	1.18658521305166\\
0.025	0.540039989937685\\
0.0186567164179104	0.0208632430971458\\
0.0140449438202247	0.00366405364688154\\
0.0105042016806723	0.00405105455949531\\
0.00788643533123028	0.000235575972283295\\
0.00592417061611374	3.6019140010415e-05\\
0.0044404973357016	1.17197186010694e-05\\
0.00333333333333333	4.38858738158706e-06\\
0.0025	1.94586103059033e-08\\
0.00187406296851574	1.22736853601332e-08\\
0.0014052838673412	5.24546624739072e-09\\
0.0010539629005059	4.55307449826848e-09\\
0.000790388871324692	4.56738393171219e-09\\
0.000592838510789661	4.54036650184906e-09\\
0.00044452347083926	4.53892534773062e-09\\
0.000333377783704494	4.53187063939086e-09\\
0.00025	4.61324556454899e-09\\
};

\addplot [color=colorclassyorange, dashed, mark=triangle,mark size=1,mark
options={solid},mark repeat=1]
  table[row sep=crcr]{%
0.125	1.41393454450555\\
0.1	1.39372137016436\\
0.0833333333333333	1.41643648417851\\
0.0714285714285714	1.41382686605232\\
0.0625	1.41327339391627\\
0.0555555555555556	1.41441213163303\\
0.05	1.41317351152251\\
0.0454545454545455	1.40922865292351\\
0.0416666666666667	1.41351822071285\\
0.0384615384615385	1.41419617603931\\
0.0357142857142857	1.41421366899153\\
0.0333333333333333	1.41420996114924\\
0.03125	1.41422073918273\\
0.0294117647058824	1.4143311162728\\
0.0277777777777778	1.41414897809977\\
0.0263157894736842	1.41421314446922\\
0.025	1.4142135504405\\
0.0186567164179104	1.41421525553721\\
0.0140449438202247	1.4193109312896\\
0.0105042016806723	0.468414368822105\\
0.00788643533123028	0.00023557596509619\\
0.00592417061611374	3.60191400029352e-05\\
0.0044404973357016	1.17197186068753e-05\\
0.00333333333333333	4.38858738447176e-06\\
0.0025	1.94586032154547e-08\\
0.00187406296851574	1.2273683368428e-08\\
0.0014052838673412	5.2454732808484e-09\\
0.0010539629005059	4.55306333391261e-09\\
0.000790388871324692	4.56737989928373e-09\\
0.000592838510789661	4.54038537471778e-09\\
0.00044452347083926	4.53892152089083e-09\\
0.000333377783704494	4.53187563510173e-09\\
0.00025	4.61324767517566e-09\\
};

\addplot [color=colorclassyorange, densely dotted, mark=square,mark
size=1,mark options={solid},mark repeat=1]
  table[row sep=crcr]{%
0.125	1.41421560334881\\
0.1	1.41418823101939\\
0.0833333333333333	1.41456133429578\\
0.0714285714285714	1.4142213083584\\
0.0625	1.41396765238398\\
0.0555555555555556	1.41422428022216\\
0.05	1.41425029949455\\
0.0454545454545455	1.41421356107075\\
0.0416666666666667	1.4142139979383\\
0.0384615384615385	1.41421321043415\\
0.0357142857142857	1.41451076543723\\
0.0333333333333333	1.43565681120671\\
0.03125	1.41549915089193\\
0.0294117647058824	1.41421281644969\\
0.0277777777777778	1.41427856029978\\
0.0263157894736842	1.41768541875128\\
0.025	1.42210610261137\\
0.0186567164179104	1.41421356237448\\
0.0140449438202247	1.41421356200633\\
0.0105042016806723	1.41421356237436\\
0.00788643533123028	1.41421356242206\\
0.00592417061611374	1.41421356237437\\
0.0044404973357016	1.4142137411191\\
0.00333333333333333	0.140054987093893\\
0.0025	3.97122125496013e-05\\
0.00187406296851574	1.2274170277536e-08\\
0.0014052838673412	5.24546696376328e-09\\
0.0010539629005059	4.55307393292485e-09\\
0.000790388871324692	4.56738450833983e-09\\
0.000592838510789661	4.54036823067494e-09\\
0.00044452347083926	4.53892647853882e-09\\
0.000333377783704494	4.53187277521389e-09\\
0.00025	4.61324571460326e-09\\
};

\end{axis}
\end{tikzpicture}%
        \hspace{-0.38cm}
%
%
\definecolor{mycolor1}{rgb}{0.00000,0.44700,0.74100}%
\definecolor{mycolor2}{rgb}{0.85000,0.32500,0.09800}%
\definecolor{mycolor3}{rgb}{0.92900,0.69400,0.12500}%
\definecolor{mycolor4}{rgb}{0.49400,0.18400,0.55600}%
\begin{tikzpicture}

\begin{axis}[%
width=1.1in, height=1.5in, scale only axis, xmode=log, xmin=0.2, xmax=400,
xminorticks=true, ymode=log, ymin=1e-9, ymax=10, yminorticks=true, axis
background/.style={fill=white}, yticklabels={}, xticklabels={}, grid=both,
title={$M=1000$}]
\addplot [color=colorclassyblue,line width=0.5pt, solid, mark=o,mark
size=1,mark options={solid},mark repeat=1]
  table[row sep=crcr]{%
0.225610947722081	1.84227775654152\\
0.267926621731598	1.17265319593959\\
0.342816668102596	0.4732466153733\\
0.400938467094363	0.27053796669742\\
0.427357057280355	0.450235287063338\\
0.493190544963537	0.386345109425166\\
0.530386707160862	0.19631810855591\\
0.567126874375033	0.0416870997865452\\
0.607334284548675	0.04468753649014\\
0.679083399050846	0.0114665093876507\\
0.716393334848313	0.0436778345414804\\
0.750803394716042	0.00214879119160973\\
0.774309500816078	0.00187304850615792\\
0.814123056205059	0.0100873067378623\\
0.925604373469197	0.000837317401088285\\
0.855817028035436	0.00583994127777149\\
0.944156374758831	0.000406289736000753\\
1.30135924929837	6.11838255075089e-05\\
1.72854159355687	6.03342676943999e-06\\
2.2690040186997	9.94019024612788e-07\\
3.01511237619616	7.91513048319724e-07\\
3.93782858130925	7.87073445620558e-07\\
5.39532723718974	7.8648980719016e-07\\
7.18786633072935	7.86692037205271e-07\\
9.58807382495737	7.86858135252198e-07\\
12.7153891852639	7.86963347129766e-07\\
17.1560310364775	7.87004767794613e-07\\
22.5176344577307	7.87039954427015e-07\\
28.7498679896025	7.87053405105163e-07\\
37.4773222797829	7.87065032174514e-07\\
49.2450133991679	7.87070991842819e-07\\
64.9983780507614	7.87074276806681e-07\\
86.3648998537147	7.87076615494501e-07\\
};

\addplot [color=black, line width=0.5pt, solid]
  table[row sep=crcr]{%
0.218386774404362	1.42142337874614\\
0.258692348079973	1.41543353882163\\
0.348904906584069	1.42618001574328\\
0.367364147867255	0.616961466424723\\
0.420587377974812	1.36502197453387\\
0.454592776066838	1.31322899363287\\
0.494744650314288	1.48065673810313\\
0.531613300644187	1.40929566131312\\
0.583880712432873	1.51482248887002\\
0.622097234701128	1.50566469199867\\
0.673078914749952	1.61449457806602\\
0.735995415434245	1.6896145376532\\
0.755784216809854	1.60202745073351\\
0.795835225807575	1.53472430191142\\
0.857043621518761	1.6715422990589\\
0.798542677257918	1.45625161453631\\
1.11260974724549	0.500329870104081\\
1.25506330044906	0.180544404447446\\
1.66682947198543	0.0199669813715889\\
2.18212931945676	0.0519286846596421\\
2.90562705013873	0.0021085765151355\\
3.81153147971422	0.00101666768566824\\
5.12680563056241	0.000323378572772979\\
6.8241183221426	8.54421422486014e-05\\
9.05348650041712	9.18676598854016e-07\\
12.1668989669417	8.34437282123676e-07\\
16.1233131712905	7.87312171250678e-07\\
21.6378502253108	7.8703713763302e-07\\
27.3666803054896	7.8705293254456e-07\\
35.7165158140801	7.87065057624439e-07\\
46.7192857059232	7.87071030196728e-07\\
61.8034346418702	7.87074294404881e-07\\
82.1568775841093	7.87076620608476e-07\\
};

\addplot [color=colorclassyorange, solid, mark=10-pointed star,mark
size=1,mark options={solid},mark repeat=1]
  table[row sep=crcr]{%
1.98612232253045	1.40415419118279\\
2.17256393160736	1.40708122746187\\
2.63113428076623	1.41475609185888\\
3.03491198456281	1.41466864823918\\
3.47353502981079	1.4143050722728\\
3.90948634671769	1.41131467275486\\
4.32928601512318	1.41438538641492\\
4.6312870796118	1.41642573597532\\
5.20401137495846	1.41323061481149\\
5.60629020913629	1.4142279935137\\
6.054084787449	1.4142122496858\\
6.39278909526573	1.40199193304616\\
6.9225096715065	1.36884212575885\\
7.20508456920782	0.713325469965504\\
7.71304436840192	1.13312138961445\\
8.11450517339667	0.18721245822371\\
8.4459012535209	0.351335182829564\\
11.4376291022663	0.0208632431077521\\
15.0840600784648	0.00366405373153227\\
20.2426814906563	0.00405105463605388\\
28.3584339709358	0.000235577286684763\\
37.5511241221946	3.60277378172235e-05\\
51.1041745905426	1.17461176932193e-05\\
68.1275807319754	4.4586064492763e-06\\
91.6974540830463	7.87305987819122e-07\\
122.319521298328	7.87161176196236e-07\\
164.615701597363	7.87082965662812e-07\\
216.165330752464	7.87078658505635e-07\\
288.431836254874	7.87078737212971e-07\\
383.075275495718	7.8707858793019e-07\\
511.379677813578	7.87078575493616e-07\\
681.962868679751	7.87078525006732e-07\\
906.026608012722	7.87079019076335e-07\\
};

\addplot [color=colorclassyorange, dashed, mark=triangle,mark size=1,mark
options={solid},mark repeat=1]
  table[row sep=crcr]{%
1.04709056288513	1.41282517662883\\
1.18059832905907	1.41443157385466\\
1.42604851078162	1.41412757202444\\
1.64886915345537	1.41096148838448\\
1.96465303404773	1.4144165082393\\
2.10868951124482	1.41693850510189\\
2.23105256790618	1.41437644821367\\
2.57320310508868	1.41418644260865\\
2.80829988919832	1.41414128624683\\
3.05522462480979	1.41295681298473\\
3.1632771239715	1.41464728915028\\
3.40849135305629	1.41423784240626\\
3.68053363720991	1.41427256940579\\
3.86702808046525	1.41421246286528\\
4.10116244198372	1.41421160656366\\
4.37644111976037	1.41421223135242\\
4.51100256756096	1.41421231095379\\
6.3478046948571	1.41416788356505\\
8.58266312477542	1.41421398941181\\
10.802368352087	0.0040513327419151\\
14.2932804231057	0.000235577286690319\\
18.5851349237822	3.60277378130524e-05\\
26.6502644560767	1.17461176915773e-05\\
35.3530385806578	4.45860644940511e-06\\
47.343788696908	7.87305987497348e-07\\
63.5279419256316	7.87161176142256e-07\\
84.2294414256582	7.87082965541765e-07\\
112.973326544942	7.87078658567248e-07\\
151.180559303814	7.87078736916059e-07\\
199.205027532311	7.87078587900957e-07\\
267.058115795107	7.87078574762062e-07\\
354.569618486199	7.87078526117217e-07\\
471.587238684705	7.8707901680075e-07\\
};

\addplot [color=colorclassyorange, densely dotted, mark=square,mark
size=1,mark options={solid},mark repeat=1]
  table[row sep=crcr]{%
0.411251038296674	1.41422465226304\\
0.543409731755473	1.41776764893643\\
0.63365861375724	1.41380764671022\\
0.733118954651764	1.41323939994499\\
0.841942952896229	1.41419217613424\\
0.930395172637179	1.41001414309911\\
1.01347571744161	1.40987510823398\\
1.11565311600074	1.41416218594551\\
1.21274792663374	1.41562183619228\\
1.27109397068964	1.41422562242581\\
1.48228448983648	1.41784190360511\\
1.47474571300174	1.41421114359345\\
1.61600209175316	1.41421356252877\\
1.65006692829591	1.41421350441019\\
1.74404122026539	1.41421347502064\\
1.87467252565685	1.41421386734806\\
1.94351666092795	1.41420604231314\\
2.49402538502157	1.41421356205934\\
3.46651378776931	1.41421356237804\\
4.63778838530645	1.41421356237033\\
6.23952224694155	1.41421356235043\\
8.28780938389044	1.4142135623013\\
11.0517624851516	1.41421566540595\\
14.5146814474284	0.000709470511975798\\
18.7107100583367	7.87254286335881e-07\\
24.4840901579232	7.87162463721254e-07\\
31.6403592483978	7.87082973911847e-07\\
41.4904593781717	7.87078659209513e-07\\
55.3778085043214	7.8707873755265e-07\\
72.6040366513096	7.87078588145297e-07\\
95.6940976117452	7.87078574197552e-07\\
126.688376210571	7.87078525410382e-07\\
168.553109306411	7.87079017918864e-07\\
};

\end{axis}
\end{tikzpicture}%
        \hspace{-0.4cm}
%
%
\definecolor{mycolor1}{rgb}{0.00000,0.44700,0.74100}%
\definecolor{mycolor2}{rgb}{0.85000,0.32500,0.09800}%
\definecolor{mycolor3}{rgb}{0.92900,0.69400,0.12500}%
\definecolor{mycolor4}{rgb}{0.49400,0.18400,0.55600}%
\begin{tikzpicture}

\begin{axis}[%
width=1.1in, height=1.5in, scale only axis, xmode=log, xmin=0.2, xmax=400,
xminorticks=true, ymode=log, ymin=1e-9, ymax=10, yminorticks=true,
yticklabels={}, axis background/.style={fill=white}, xticklabels={},
grid=both, title={$M=2000$}]
\addplot [color=colorclassyblue,line width=0.5pt, solid, mark=o,mark
size=1,mark options={solid},mark repeat=1]
  table[row sep=crcr]{%
0.289765449269132	1.85494302749674\\
0.386437886668886	1.35892792944903\\
0.453549301042845	0.488130494494647\\
0.499809826977071	0.26497196668039\\
0.5722274739529	0.459677280081579\\
0.632701594661587	0.380603252191244\\
0.698976968200729	0.198323808337455\\
0.749782734450916	0.0524171094127853\\
0.803880130056133	0.0447008597341876\\
0.880151069727006	0.0114680381560209\\
0.949412174735837	0.0436777561208933\\
0.97650650205618	0.00214763527462076\\
1.1361728902427	0.00187279773846845\\
1.02522771981195	0.0100874432412646\\
1.10843914928567	0.000837311912421154\\
1.22569500209686	0.00583994129926332\\
1.18212541693221	0.000406288999606731\\
1.59832425596118	6.11789073697043e-05\\
2.10212996787622	5.98275056956967e-06\\
2.81311530429035	6.11907841941008e-07\\
3.99657268286855	1.02235240051469e-07\\
5.46397804739886	4.37777364948667e-08\\
7.06985608990456	1.49455166558599e-08\\
8.7929380516256	7.40374031154423e-09\\
12.4906584069235	5.30796338809436e-09\\
16.0229192264476	4.76959759369982e-09\\
21.8545430038692	4.63303159501684e-09\\
29.2663641109433	4.58366429543514e-09\\
36.6605911364062	4.57348720372152e-09\\
51.6398816394326	4.54248469051601e-09\\
67.0670225568921	4.53968159488868e-09\\
89.3631052928142	4.53222928927586e-09\\
119.978461366658	4.61340414162414e-09\\
};

\addplot [color=black, line width=0.5pt, solid]
  table[row sep=crcr]{%
0.27438139732593	1.4193637600574\\
0.375778891364824	1.51409377932297\\
0.430757056740005	1.44034071060893\\
0.554608476223269	1.36201421285788\\
0.535169701279758	1.14941806946162\\
0.599850743445935	0.941130481007499\\
0.649170548427778	1.2270336353767\\
0.718909862790239	1.26519672881492\\
0.774238955180106	1.33570829005342\\
0.825005996379058	1.41713762884927\\
0.877080383299957	0.938311311543404\\
0.942915071759922	1.77148898716842\\
1.09131457139624	1.6166610234424\\
1.01245325601131	1.53954904178951\\
1.10837130539746	1.64889286193774\\
1.16217961399229	1.47051166807726\\
1.23767395127908	0.502688848888131\\
1.65243095758646	0.1767862286821\\
2.17671831908061	0.0199669813553296\\
2.88591930119594	0.0519286846538669\\
3.75817436299542	0.00210857636874647\\
5.11025622475179	0.00101666738190946\\
6.69811090794635	0.000323377616812039\\
8.98501760789055	8.54385208755247e-05\\
12.0172881837261	4.74182934984492e-07\\
16.08929786679	2.77483536405477e-07\\
21.2792595530049	2.24814075849563e-08\\
27.3241956221479	4.07126373065946e-09\\
35.3904109568478	4.49145664427237e-09\\
46.8656922162947	4.5468056509284e-09\\
61.6319559122793	4.54614140448793e-09\\
81.8939858194264	4.53513857704442e-09\\
108.811054471136	4.61449291227918e-09\\
};

\addplot [color=colorclassyorange, solid, mark=10-pointed star,mark
size=1,mark options={solid},mark repeat=1]
  table[row sep=crcr]{%
2.57650974415349	1.4295074309552\\
3.13481901442043	1.40729746572539\\
3.7669244234095	1.41721791013935\\
4.26598736722786	1.41240231872207\\
4.96175525947736	1.41526693378122\\
5.48249192402565	1.4142365838197\\
6.28126875275616	1.41439238165919\\
6.65028276792257	1.41419980934031\\
7.30898478225564	1.41427082254659\\
7.80098715859289	1.41421400623906\\
8.65958128313206	1.4142201367728\\
8.98231766129509	1.41418020835033\\
9.69810100153789	1.41421348023915\\
10.3429325430621	1.41422876944897\\
10.8491560190287	1.40992447977729\\
11.3722918355884	1.18658521305166\\
12.3438334258432	0.540039989937685\\
16.2399749397885	0.0208632430971458\\
22.0767947785422	0.00366405364688154\\
29.8183987221333	0.00405105455949531\\
39.7271835150159	0.000235575972283295\\
53.5345200810404	3.6019140010415e-05\\
71.3219249532522	1.17197186010694e-05\\
94.9195785753804	4.38858738158706e-06\\
127.100798306424	1.94586103059033e-08\\
169.521438218384	1.22736853601332e-08\\
227.431741693702	5.24546624739072e-09\\
301.032387951766	4.55307449826848e-09\\
401.345652902863	4.56738393171219e-09\\
536.483354682144	4.54036650184906e-09\\
711.430332882337	4.53892534773062e-09\\
955.601760578921	4.53187063939086e-09\\
1311.75962550173	4.61324556454899e-09\\
};

\addplot [color=colorclassyorange, dashed, mark=triangle,mark size=1,mark
options={solid},mark repeat=1]
  table[row sep=crcr]{%
1.40777869203588	1.41393454450555\\
1.66284859675724	1.39372137016436\\
2.00975931329977	1.41643648417851\\
2.30861974611376	1.41382686605232\\
2.57905749229176	1.41327339391627\\
2.81711929447159	1.41441213163303\\
3.242743931049	1.41317351152251\\
3.55667411248337	1.40922865292351\\
3.81573029592844	1.41351822071285\\
4.19153650496506	1.41419617603931\\
4.51921017706354	1.41421366899153\\
4.74620021692032	1.41420996114924\\
5.1707285442953	1.41422073918273\\
5.41595538153632	1.4143311162728\\
5.7724140746649	1.41414897809977\\
6.07634298627782	1.41421314446922\\
6.29367637924974	1.4142135504405\\
8.69734000925199	1.41421525553721\\
11.4954558103637	1.4193109312896\\
15.2387057187294	0.468414368822105\\
20.6126450200635	0.00023557596509619\\
27.0577900844056	3.60191400029352e-05\\
37.5416614007601	1.17197186068753e-05\\
49.9810907676176	4.38858738447176e-06\\
67.1037792347871	1.94586032154547e-08\\
89.5216897811194	1.2273683368428e-08\\
120.945039546082	5.2454732808484e-09\\
159.172730239192	4.55306333391261e-09\\
212.416599898475	4.56737989928373e-09\\
283.318356845043	4.54038537471778e-09\\
375.388963120242	4.53892152089083e-09\\
499.391357469148	4.53187563510173e-09\\
667.602455828676	4.61324767517566e-09\\
};

\addplot [color=colorclassyorange, densely dotted, mark=square,mark
size=1,mark options={solid},mark repeat=1]
  table[row sep=crcr]{%
0.657482925704638	1.41421560334881\\
0.767474078981693	1.41418823101939\\
0.92762017749282	1.41456133429578\\
1.0648236344091	1.4142213083584\\
1.2170392026602	1.41396765238398\\
1.35188133203369	1.41422428022216\\
1.51018844090242	1.41425029949455\\
1.62281710038224	1.41421356107075\\
1.8561214249498	1.4142139979383\\
1.83359605328687	1.41421321043415\\
2.07858693789012	1.41451076543723\\
2.19534566911335	1.43565681120671\\
2.34232675117534	1.41549915089193\\
2.4906475999324	1.41421281644969\\
2.52674535156793	1.41427856029978\\
2.78033019560354	1.41768541875128\\
2.93500075198646	1.42210610261137\\
3.96794586417954	1.41421356237448\\
5.10519825274971	1.41421356200633\\
6.98929147265366	1.41421356237436\\
9.05648093751848	1.41421356242206\\
12.3231527475124	1.41421356237437\\
16.3145692948769	1.4142137411191\\
21.7580942111452	0.140054987093893\\
28.4151781787638	3.97122125496013e-05\\
39.1944754061402	1.2274170277536e-08\\
52.3257863511906	5.24546696376328e-09\\
69.1760750029194	4.55307393292485e-09\\
92.7917159610557	4.56738450833983e-09\\
123.694412515456	4.54036823067494e-09\\
163.756366293355	4.53892647853882e-09\\
217.174728031562	4.53187277521389e-09\\
289.519253405027	4.61324571460326e-09\\
};

\end{axis}
\end{tikzpicture}%
\\
        \vspace{-0.2cm}
%
%
\definecolor{mycolor1}{rgb}{0.00000,0.44700,0.74100}%
\definecolor{mycolor2}{rgb}{0.85000,0.32500,0.09800}%
\definecolor{mycolor3}{rgb}{0.92900,0.69400,0.12500}%
\definecolor{mycolor4}{rgb}{0.49400,0.18400,0.55600}%
\begin{tikzpicture}

\begin{axis}[%
width=1.1in, height=1.5in, scale only axis, xmode=log, xmin=0.0005, xmax=0.2,
xminorticks=true, ymode=log, ymin=1e-12, ymax=50, yminorticks=true, axis
background/.style={fill=white}, xticklabel style =
{xshift=-0.25ex,yshift=-0.1ex}, ylabel={Error in Total Energy}, xlabel={time
step \mbox{}}, grid=both]
\addplot [color=colorclassyblue,line width=0.5pt, solid, mark=o,mark
size=1,mark options={solid},mark repeat=1]
  table[row sep=crcr]{%
0.125	0.169741794901054\\
0.1	0.0204417836300359\\
0.0833333333333333	0.024235097798865\\
0.0714285714285714	0.0629817107622301\\
0.0625	0.0291402313010143\\
0.0555555555555556	0.0704055895671418\\
0.05	0.0618998546768923\\
0.0454545454545455	0.00947477377839157\\
0.0416666666666667	0.000436456111272854\\
0.0384615384615385	0.000117286470124434\\
0.0357142857142857	0.000461381612276934\\
0.0333333333333333	2.36565226918728e-05\\
0.03125	3.11332208520199e-05\\
0.0294117647058824	0.000122000837898817\\
0.0277777777777778	1.53845070727243e-06\\
0.0263157894736842	7.06214924033333e-05\\
0.025	1.31635729472279e-06\\
0.0186567164179104	7.44901011540122e-07\\
0.0140449438202247	7.70246439962818e-07\\
0.0105042016806723	1.12608649072143e-07\\
0.00788643533123028	5.32637929318014e-08\\
0.00592417061611374	1.33607933605617e-08\\
0.0044404973357016	4.16852730111827e-09\\
0.00333333333333333	1.14337872503256e-09\\
0.0025	1.30190969116484e-10\\
0.00187406296851574	3.18739257210154e-10\\
0.0014052838673412	6.66222632617064e-12\\
0.0010539629005059	3.00510283324229e-10\\
0.000790388871324692	4.01361610613549e-10\\
0.000592838510789661	2.17542428515571e-10\\
0.00044452347083926	8.34334379362645e-10\\
0.000333377783704494	1.12286802078643e-09\\
0.00025	1.14703357922963e-09\\
};

\addplot [color=black, line width=0.5pt, solid]
  table[row sep=crcr]{%
0.125	2.76884638360543\\
0.1	5.85830325188471\\
0.0833333333333333	1.29133854367552\\
0.0714285714285714	1.73660172615262\\
0.0625	17.8793592876427\\
0.0555555555555556	1.881013385979\\
0.05	0.878583516734298\\
0.0454545454545455	7.36053318088371\\
0.0416666666666667	25.6286235331592\\
0.0384615384615385	43.6181911292148\\
0.0357142857142857	23.7783739754883\\
0.0333333333333333	7.65201767637111\\
0.03125	2.34356158901083\\
0.0294117647058824	1.82289893568458\\
0.0277777777777778	1.3465399587208\\
0.0263157894736842	2.00152071220993\\
0.025	0.712383693314172\\
0.0186567164179104	0.163917932808825\\
0.0140449438202247	0.000668852473412862\\
0.0105042016806723	0.000541067913569115\\
0.00788643533123028	3.91211232830813e-06\\
0.00592417061611374	7.46575372900793e-06\\
0.0044404973357016	4.05386459512158e-06\\
0.00333333333333333	1.10795993890633e-06\\
0.0025	3.45126061063183e-09\\
0.00187406296851574	5.04326935768518e-09\\
0.0014052838673412	2.2336825722391e-09\\
0.0010539629005059	1.86543580582565e-09\\
0.000790388871324692	1.33526967260877e-10\\
0.000592838510789661	1.02542685453955e-09\\
0.00044452347083926	1.48268730626455e-09\\
0.000333377783704494	1.51113077606624e-09\\
0.00025	1.35292532732478e-09\\
};

\addplot [color=colorclassyorange, solid, mark=10-pointed star,mark
size=1,mark options={solid},mark repeat=1]
  table[row sep=crcr]{%
0.125	7.55770506082554\\
0.1	4.3590425301067\\
0.0833333333333333	11.6855835873345\\
0.0714285714285714	13.5267579193662\\
0.0625	9.25328147396188\\
0.0555555555555556	13.2302944241713\\
0.05	10.1829504928097\\
0.0454545454545455	13.1018393316185\\
0.0416666666666667	7.36784686972454\\
0.0384615384615385	7.29226234321017\\
0.0357142857142857	1.29346449623165\\
0.0333333333333333	11.2326797807017\\
0.03125	6.86580153867713\\
0.0294117647058824	1.98867085470199\\
0.0277777777777778	0.0728330037157985\\
0.0263157894736842	0.00343889828734056\\
0.025	0.00320906485425132\\
0.0186567164179104	0.000152399784570711\\
0.0140449438202247	2.46411885242992e-05\\
0.0105042016806723	3.59437191055356e-05\\
0.00788643533123028	1.88838538939251e-06\\
0.00592417061611374	7.08343024236058e-08\\
0.0044404973357016	1.31552778448452e-07\\
0.00333333333333333	5.8717942508224e-08\\
0.0025	1.25457555455455e-09\\
0.00187406296851574	3.17374126979075e-10\\
0.0014052838673412	2.98792102171319e-11\\
0.0010539629005059	3.11196401980851e-10\\
0.000790388871324692	4.06092937055291e-10\\
0.000592838510789661	2.16037854272599e-10\\
0.00044452347083926	8.33838775804452e-10\\
0.000333377783704494	1.1226819474075e-09\\
0.00025	1.14693365915741e-09\\
};

\addplot [color=colorclassyorange, dashed, mark=triangle,mark size=1,mark
options={solid},mark repeat=1]
  table[row sep=crcr]{%
0.125	6.93074680842304\\
0.1	7.58342555040532\\
0.0833333333333333	11.6477603772651\\
0.0714285714285714	4.15651766335908\\
0.0625	10.1182997517946\\
0.0555555555555556	12.4048789049898\\
0.05	14.2414402384593\\
0.0454545454545455	15.7659772262573\\
0.0416666666666667	19.9549136622432\\
0.0384615384615385	7.65034055316803\\
0.0357142857142857	10.1989036504273\\
0.0333333333333333	4.7322366944302\\
0.03125	9.99386202712216\\
0.0294117647058824	15.8269919724448\\
0.0277777777777778	7.69375992700619\\
0.0263157894736842	8.67288126757482\\
0.025	19.5043768910508\\
0.0186567164179104	15.4560800118047\\
0.0140449438202247	4.30398667918028\\
0.0105042016806723	3.59803073637366e-05\\
0.00788643533123028	1.88838538761615e-06\\
0.00592417061611374	7.08343019795166e-08\\
0.0044404973357016	1.31552778892541e-07\\
0.00333333333333333	5.87179433964025e-08\\
0.0025	1.25457599864376e-09\\
0.00187406296851574	3.17373682889865e-10\\
0.0014052838673412	2.98792102171319e-11\\
0.0010539629005059	3.11197734248481e-10\\
0.000790388871324692	4.06093825233711e-10\\
0.000592838510789661	2.16038298361809e-10\\
0.00044452347083926	8.33839219893662e-10\\
0.000333377783704494	1.12268372376434e-09\\
0.00025	1.14693188280057e-09\\
};

\addplot [color=colorclassyorange, densely dotted, mark=square,mark
size=1,mark options={solid},mark repeat=1]
  table[row sep=crcr]{%
0.125	22.5260206380454\\
0.1	18.4660317308899\\
0.0833333333333333	17.0755267562189\\
0.0714285714285714	7.60325159137131\\
0.0625	5.01505845565192\\
0.0555555555555556	14.8254805727019\\
0.05	22.7144434603609\\
0.0454545454545455	20.4009284867757\\
0.0416666666666667	21.1026226905981\\
0.0384615384615385	16.8102943373415\\
0.0357142857142857	18.2182955951569\\
0.0333333333333333	28.6261940072811\\
0.03125	23.0443258440159\\
0.0294117647058824	15.3960743712316\\
0.0277777777777778	18.9496751779985\\
0.0263157894736842	18.0312803112173\\
0.025	18.7199105988621\\
0.0186567164179104	21.6705019827711\\
0.0140449438202247	21.7687443813633\\
0.0105042016806723	17.7105290349124\\
0.00788643533123028	18.254813810078\\
0.00592417061611374	11.8848033366774\\
0.0044404973357016	3.70136574855487\\
0.00333333333333333	2.150941658563e-06\\
0.0025	1.24747900898114e-09\\
0.00187406296851574	3.17377679692754e-10\\
0.0014052838673412	2.98809865739713e-11\\
0.0010539629005059	3.11195957891641e-10\\
0.000790388871324692	4.06094269322921e-10\\
0.000592838510789661	2.16037854272599e-10\\
0.00044452347083926	8.33839219893662e-10\\
0.000333377783704494	1.12268461194276e-09\\
0.00025	1.14693410324662e-09\\
};

\end{axis}
\end{tikzpicture}%
        \hspace{-0.46cm}
%
%
\definecolor{mycolor1}{rgb}{0.00000,0.44700,0.74100}%
\definecolor{mycolor2}{rgb}{0.85000,0.32500,0.09800}%
\definecolor{mycolor3}{rgb}{0.92900,0.69400,0.12500}%
\definecolor{mycolor4}{rgb}{0.49400,0.18400,0.55600}%
\begin{tikzpicture}

\begin{axis}[%
width=1.1in, height=1.5in, scale only axis, xmode=log, xmin=0.0005, xmax=0.2,
xminorticks=true, ymode=log, ymin=1e-12, ymax=50, yticklabels={},
yminorticks=true, axis background/.style={fill=white}, xticklabel style =
{xshift=-0.25ex, yshift=-0.1ex}, xlabel={time step \mbox{}},  grid=both]
\addplot [color=colorclassyblue,line width=0.5pt, solid, mark=o,mark
size=1,mark options={solid},mark repeat=1]
  table[row sep=crcr]{%
0.125	0.225069627745611\\
0.1	0.123263546797258\\
0.0833333333333333	0.0246505898034428\\
0.0714285714285714	0.0610825727742919\\
0.0625	0.0641589663215911\\
0.0555555555555556	0.130657881357833\\
0.05	0.0802444676283434\\
0.0454545454545455	0.0201127655087898\\
0.0416666666666667	0.000501937216514037\\
0.0384615384615385	0.000112814058672583\\
0.0357142857142857	0.000459603809048303\\
0.0333333333333333	1.93949291955597e-05\\
0.03125	2.92879960177572e-05\\
0.0294117647058824	0.00012224968162311\\
0.0277777777777778	1.52898537653101e-06\\
0.0263157894736842	7.06217217611993e-05\\
0.025	1.31641350264999e-06\\
0.0186567164179104	7.44899150806333e-07\\
0.0140449438202247	7.70246330716873e-07\\
0.0105042016806723	1.12608651292589e-07\\
0.00788643533123028	5.32638018135856e-08\\
0.00592417061611374	1.33608026864351e-08\\
0.0044404973357016	4.16853485063484e-09\\
0.00333333333333333	1.14338583045992e-09\\
0.0025	1.302007390791e-10\\
0.00187406296851574	3.18729043158328e-10\\
0.0014052838673412	6.65245636355394e-12\\
0.0010539629005059	3.00498737004773e-10\\
0.000790388871324692	4.01348287937253e-10\\
0.000592838510789661	2.17558415727126e-10\\
0.00044452347083926	8.34357027912347e-10\\
0.000333377783704494	1.12289200160376e-09\\
0.00025	1.147068218188e-09\\
};

\addplot [color=black, line width=0.5pt, solid]
  table[row sep=crcr]{%
0.125	69.143024572697\\
0.1	68.6930680566039\\
0.0833333333333333	43.2983841758332\\
0.0714285714285714	79.7825938878128\\
0.0625	77.1434023842689\\
0.0555555555555556	41.4190581387187\\
0.05	59.0316036272882\\
0.0454545454545455	49.489979966796\\
0.0416666666666667	47.9719347552226\\
0.0384615384615385	48.3855438127671\\
0.0357142857142857	40.9006432667609\\
0.0333333333333333	3.45121625246258\\
0.03125	1.13245703237228\\
0.0294117647058824	3.22040654995179\\
0.0277777777777778	2.52197411969148\\
0.0263157894736842	2.21352813735734\\
0.025	0.806177018749089\\
0.0186567164179104	0.114260779897472\\
0.0140449438202247	0.000668852473520776\\
0.0105042016806723	0.000541067913612636\\
0.00788643533123028	3.91211231898225e-06\\
0.00592417061611374	7.46575374011016e-06\\
0.0044404973357016	4.05386460400337e-06\\
0.00333333333333333	1.10795994601176e-06\\
0.0025	3.45126860423761e-09\\
0.00187406296851574	5.04326047590098e-09\\
0.0014052838673412	2.23367235818728e-09\\
0.0010539629005059	1.86542381541699e-09\\
0.000790388871324692	1.33511868227743e-10\\
0.000592838510789661	1.02544417401873e-09\\
0.00044452347083926	1.48270862254662e-09\\
0.000333377783704494	1.51115919777567e-09\\
0.00025	1.3529577458371e-09\\
};

\addplot [color=colorclassyorange, solid, mark=10-pointed star,mark
size=1,mark options={solid},mark repeat=1]
  table[row sep=crcr]{%
0.125	10.0293607498953\\
0.1	6.20552505122641\\
0.0833333333333333	9.25487900534435\\
0.0714285714285714	9.29507189108535\\
0.0625	8.31274025079756\\
0.0555555555555556	11.5429362858658\\
0.05	14.645432835166\\
0.0454545454545455	12.6869721388933\\
0.0416666666666667	8.57584186661356\\
0.0384615384615385	11.328605484234\\
0.0357142857142857	5.16659098328052\\
0.0333333333333333	10.4170433023926\\
0.03125	1.95100496292573\\
0.0294117647058824	1.95757567510364\\
0.0277777777777778	4.35533925918163\\
0.0263157894736842	4.35746803924623\\
0.025	0.503098588137472\\
0.0186567164179104	0.000152399853704743\\
0.0140449438202247	2.46411885198583e-05\\
0.0105042016806723	3.59437191037593e-05\\
0.00788643533123028	1.88838538228708e-06\\
0.00592417061611374	7.08342984268029e-08\\
0.0044404973357016	1.31552774895738e-07\\
0.00333333333333333	5.87179385114212e-08\\
0.0025	1.25457511046534e-09\\
0.00187406296851574	3.17371906533026e-10\\
0.0014052838673412	2.98747693250334e-11\\
0.0010539629005059	3.11195069713222e-10\\
0.000790388871324692	4.06092048876872e-10\\
0.000592838510789661	2.16041406986278e-10\\
0.00044452347083926	8.3384410487497e-10\\
0.000333377783704494	1.12268150331829e-09\\
0.00025	1.14693721187109e-09\\
};

\addplot [color=colorclassyorange, dashed, mark=triangle,mark size=1,mark
options={solid},mark repeat=1]
  table[row sep=crcr]{%
0.125	10.9664944266125\\
0.1	14.3550275155756\\
0.0833333333333333	11.6107201964703\\
0.0714285714285714	14.8852474187719\\
0.0625	12.3961577836211\\
0.0555555555555556	17.0629311837632\\
0.05	12.7014888837083\\
0.0454545454545455	21.2932110228841\\
0.0416666666666667	18.1097839275809\\
0.0384615384615385	20.9260509945046\\
0.0357142857142857	17.1018659089637\\
0.0333333333333333	18.7766130937429\\
0.03125	17.1486838650711\\
0.0294117647058824	16.3529525435747\\
0.0277777777777778	17.395277791378\\
0.0263157894736842	18.041755304413\\
0.025	21.359779347262\\
0.0186567164179104	19.6541496426341\\
0.0140449438202247	18.3614209023073\\
0.0105042016806723	0.327511886318357\\
0.00788643533123028	1.88838300418936e-06\\
0.00592417061611374	7.08342984268029e-08\\
0.0044404973357016	1.31552774451649e-07\\
0.00333333333333333	5.87179407318672e-08\\
0.0025	1.25457377819771e-09\\
0.00187406296851574	3.17372794711446e-10\\
0.0014052838673412	2.98752134142433e-11\\
0.0010539629005059	3.11194181534802e-10\\
0.000790388871324692	4.06092492966081e-10\\
0.000592838510789661	2.16041406986278e-10\\
0.00044452347083926	8.33842772607341e-10\\
0.000333377783704494	1.12268061513987e-09\\
0.00025	1.1469376559603e-09\\
};

\addplot [color=colorclassyorange, densely dotted, mark=square,mark
size=1,mark options={solid},mark repeat=1]
  table[row sep=crcr]{%
0.125	1.10720383783545\\
0.1	16.9045006987546\\
0.0833333333333333	15.601588622683\\
0.0714285714285714	24.5583961411377\\
0.0625	16.4003066629779\\
0.0555555555555556	25.5672555682025\\
0.05	21.1069795121477\\
0.0454545454545455	26.6706129707234\\
0.0416666666666667	22.7901674174583\\
0.0384615384615385	24.4796721974552\\
0.0357142857142857	24.6195358045807\\
0.0333333333333333	20.307999337318\\
0.03125	24.4373986516784\\
0.0294117647058824	23.2772585242079\\
0.0277777777777778	20.2325084090028\\
0.0263157894736842	29.0279988657505\\
0.025	29.387300246597\\
0.0186567164179104	12.6165719901597\\
0.0140449438202247	27.7657139669748\\
0.0105042016806723	2.34689623148446\\
0.00788643533123028	18.9626453575675\\
0.00592417061611374	20.4786675537396\\
0.0044404973357016	9.89931983490108\\
0.00333333333333333	0.0207526901844393\\
0.0025	6.68552546656542e-09\\
0.00187406296851574	3.16660031529636e-10\\
0.0014052838673412	2.98778779495024e-11\\
0.0010539629005059	3.11196401980851e-10\\
0.000790388871324692	4.06090716609242e-10\\
0.000592838510789661	2.16041406986278e-10\\
0.00044452347083926	8.3384366078576e-10\\
0.000333377783704494	1.1226819474075e-09\\
0.00025	1.1469376559603e-09\\
};

\end{axis}
\end{tikzpicture}%
        \hspace{-0.46cm}
%
%
\definecolor{mycolor1}{rgb}{0.00000,0.44700,0.74100}%
\definecolor{mycolor2}{rgb}{0.85000,0.32500,0.09800}%
\definecolor{mycolor3}{rgb}{0.92900,0.69400,0.12500}%
\definecolor{mycolor4}{rgb}{0.49400,0.18400,0.55600}%
\begin{tikzpicture}

\begin{axis}[%
width=1.1in, height=1.5in, scale only axis, xmode=log, xmin=0.2, xmax=400,
xminorticks=true, ymode=log, ymin=1e-12, ymax=50, yminorticks=true, axis
background/.style={fill=white}, yticklabels={}, xticklabel style =
{xshift=0ex, text height=2ex}, xlabel={cost in seconds \mbox{}}, grid=both]
\addplot [color=colorclassyblue,line width=0.5pt, solid, mark=o,mark
size=1,mark options={solid},mark repeat=1]
  table[row sep=crcr]{%
0.225610947722081	0.169741794901054\\
0.267926621731598	0.0204417836300359\\
0.342816668102596	0.024235097798865\\
0.400938467094363	0.0629817107622301\\
0.427357057280355	0.0291402313010143\\
0.493190544963537	0.0704055895671418\\
0.530386707160862	0.0618998546768923\\
0.567126874375033	0.00947477377839157\\
0.607334284548675	0.000436456111272854\\
0.679083399050846	0.000117286470124434\\
0.716393334848313	0.000461381612276934\\
0.750803394716042	2.36565226918728e-05\\
0.774309500816078	3.11332208520199e-05\\
0.814123056205059	0.000122000837898817\\
0.925604373469197	1.53845070727243e-06\\
0.855817028035436	7.06214924033333e-05\\
0.944156374758831	1.31635729472279e-06\\
1.30135924929837	7.44901011540122e-07\\
1.72854159355687	7.70246439962818e-07\\
2.2690040186997	1.12608649072143e-07\\
3.01511237619616	5.32637929318014e-08\\
3.93782858130925	1.33607933605617e-08\\
5.39532723718974	4.16852730111827e-09\\
7.18786633072935	1.14337872503256e-09\\
9.58807382495737	1.30190969116484e-10\\
12.7153891852639	3.18739257210154e-10\\
17.1560310364775	6.66222632617064e-12\\
22.5176344577307	3.00510283324229e-10\\
28.7498679896025	4.01361610613549e-10\\
37.4773222797829	2.17542428515571e-10\\
49.2450133991679	8.34334379362645e-10\\
64.9983780507614	1.12286802078643e-09\\
86.3648998537147	1.14703357922963e-09\\
};

\addplot [color=black, line width=0.5pt, solid]
  table[row sep=crcr]{%
0.218386774404362	2.76884638360543\\
0.258692348079973	5.85830325188471\\
0.348904906584069	1.29133854367552\\
0.367364147867255	1.73660172615262\\
0.420587377974812	17.8793592876427\\
0.454592776066838	1.881013385979\\
0.494744650314288	0.878583516734298\\
0.531613300644187	7.36053318088371\\
0.583880712432873	25.6286235331592\\
0.622097234701128	43.6181911292148\\
0.673078914749952	23.7783739754883\\
0.735995415434245	7.65201767637111\\
0.755784216809854	2.34356158901083\\
0.795835225807575	1.82289893568458\\
0.857043621518761	1.3465399587208\\
0.798542677257918	2.00152071220993\\
1.11260974724549	0.712383693314172\\
1.25506330044906	0.163917932808825\\
1.66682947198543	0.000668852473412862\\
2.18212931945676	0.000541067913569115\\
2.90562705013873	3.91211232830813e-06\\
3.81153147971422	7.46575372900793e-06\\
5.12680563056241	4.05386459512158e-06\\
6.8241183221426	1.10795993890633e-06\\
9.05348650041712	3.45126061063183e-09\\
12.1668989669417	5.04326935768518e-09\\
16.1233131712905	2.2336825722391e-09\\
21.6378502253108	1.86543580582565e-09\\
27.3666803054896	1.33526967260877e-10\\
35.7165158140801	1.02542685453955e-09\\
46.7192857059232	1.48268730626455e-09\\
61.8034346418702	1.51113077606624e-09\\
82.1568775841093	1.35292532732478e-09\\
};

\addplot [color=colorclassyorange, solid, mark=10-pointed star,mark
size=1,mark options={solid},mark repeat=1]
  table[row sep=crcr]{%
1.98612232253045	7.55770506082554\\
2.17256393160736	4.3590425301067\\
2.63113428076623	11.6855835873345\\
3.03491198456281	13.5267579193662\\
3.47353502981079	9.25328147396188\\
3.90948634671769	13.2302944241713\\
4.32928601512318	10.1829504928097\\
4.6312870796118	13.1018393316185\\
5.20401137495846	7.36784686972454\\
5.60629020913629	7.29226234321017\\
6.054084787449	1.29346449623165\\
6.39278909526573	11.2326797807017\\
6.9225096715065	6.86580153867713\\
7.20508456920782	1.98867085470199\\
7.71304436840192	0.0728330037157985\\
8.11450517339667	0.00343889828734056\\
8.4459012535209	0.00320906485425132\\
11.4376291022663	0.000152399784570711\\
15.0840600784648	2.46411885242992e-05\\
20.2426814906563	3.59437191055356e-05\\
28.3584339709358	1.88838538939251e-06\\
37.5511241221946	7.08343024236058e-08\\
51.1041745905426	1.31552778448452e-07\\
68.1275807319754	5.8717942508224e-08\\
91.6974540830463	1.25457555455455e-09\\
122.319521298328	3.17374126979075e-10\\
164.615701597363	2.98792102171319e-11\\
216.165330752464	3.11196401980851e-10\\
288.431836254874	4.06092937055291e-10\\
383.075275495718	2.16037854272599e-10\\
511.379677813578	8.33838775804452e-10\\
681.962868679751	1.1226819474075e-09\\
906.026608012722	1.14693365915741e-09\\
};

\addplot [color=colorclassyorange, dashed, mark=triangle,mark size=1,mark
options={solid},mark repeat=1]
  table[row sep=crcr]{%
1.04709056288513	6.93074680842304\\
1.18059832905907	7.58342555040532\\
1.42604851078162	11.6477603772651\\
1.64886915345537	4.15651766335908\\
1.96465303404773	10.1182997517946\\
2.10868951124482	12.4048789049898\\
2.23105256790618	14.2414402384593\\
2.57320310508868	15.7659772262573\\
2.80829988919832	19.9549136622432\\
3.05522462480979	7.65034055316803\\
3.1632771239715	10.1989036504273\\
3.40849135305629	4.7322366944302\\
3.68053363720991	9.99386202712216\\
3.86702808046525	15.8269919724448\\
4.10116244198372	7.69375992700619\\
4.37644111976037	8.67288126757482\\
4.51100256756096	19.5043768910508\\
6.3478046948571	15.4560800118047\\
8.58266312477542	4.30398667918028\\
10.802368352087	3.59803073637366e-05\\
14.2932804231057	1.88838538761615e-06\\
18.5851349237822	7.08343019795166e-08\\
26.6502644560767	1.31552778892541e-07\\
35.3530385806578	5.87179433964025e-08\\
47.343788696908	1.25457599864376e-09\\
63.5279419256316	3.17373682889865e-10\\
84.2294414256582	2.98792102171319e-11\\
112.973326544942	3.11197734248481e-10\\
151.180559303814	4.06093825233711e-10\\
199.205027532311	2.16038298361809e-10\\
267.058115795107	8.33839219893662e-10\\
354.569618486199	1.12268372376434e-09\\
471.587238684705	1.14693188280057e-09\\
};

\addplot [color=colorclassyorange, densely dotted, mark=square,mark
size=1,mark options={solid},mark repeat=1]
  table[row sep=crcr]{%
0.411251038296674	22.5260206380454\\
0.543409731755473	18.4660317308899\\
0.63365861375724	17.0755267562189\\
0.733118954651764	7.60325159137131\\
0.841942952896229	5.01505845565192\\
0.930395172637179	14.8254805727019\\
1.01347571744161	22.7144434603609\\
1.11565311600074	20.4009284867757\\
1.21274792663374	21.1026226905981\\
1.27109397068964	16.8102943373415\\
1.48228448983648	18.2182955951569\\
1.47474571300174	28.6261940072811\\
1.61600209175316	23.0443258440159\\
1.65006692829591	15.3960743712316\\
1.74404122026539	18.9496751779985\\
1.87467252565685	18.0312803112173\\
1.94351666092795	18.7199105988621\\
2.49402538502157	21.6705019827711\\
3.46651378776931	21.7687443813633\\
4.63778838530645	17.7105290349124\\
6.23952224694155	18.254813810078\\
8.28780938389044	11.8848033366774\\
11.0517624851516	3.70136574855487\\
14.5146814474284	2.150941658563e-06\\
18.7107100583367	1.24747900898114e-09\\
24.4840901579232	3.17377679692754e-10\\
31.6403592483978	2.98809865739713e-11\\
41.4904593781717	3.11195957891641e-10\\
55.3778085043214	4.06094269322921e-10\\
72.6040366513096	2.16037854272599e-10\\
95.6940976117452	8.33839219893662e-10\\
126.688376210571	1.12268461194276e-09\\
168.553109306411	1.14693410324662e-09\\
};

\end{axis}
\end{tikzpicture}%
        \hspace{-0.4cm}
%
%
\definecolor{mycolor1}{rgb}{0.00000,0.44700,0.74100}%
\definecolor{mycolor2}{rgb}{0.85000,0.32500,0.09800}%
\definecolor{mycolor3}{rgb}{0.92900,0.69400,0.12500}%
\definecolor{mycolor4}{rgb}{0.49400,0.18400,0.55600}%
\begin{tikzpicture}

\begin{axis}[%
width=1.1in, height=1.5in, scale only axis, xmode=log, xmin=0.2, xmax=400,
xminorticks=true, ymode=log, ymin=1e-12, ymax=50, yminorticks=true,
yticklabels={}, axis background/.style={fill=white}, xticklabel style =
{xshift=0ex, text height=2ex}, xlabel={cost in seconds \mbox{}}, grid=both]
\addplot [color=colorclassyblue,line width=0.5pt, solid, mark=o,mark
size=1,mark options={solid},mark repeat=1]
  table[row sep=crcr]{%
0.289765449269132	0.225069627745611\\
0.386437886668886	0.123263546797258\\
0.453549301042845	0.0246505898034428\\
0.499809826977071	0.0610825727742919\\
0.5722274739529	0.0641589663215911\\
0.632701594661587	0.130657881357833\\
0.698976968200729	0.0802444676283434\\
0.749782734450916	0.0201127655087898\\
0.803880130056133	0.000501937216514037\\
0.880151069727006	0.000112814058672583\\
0.949412174735837	0.000459603809048303\\
0.97650650205618	1.93949291955597e-05\\
1.1361728902427	2.92879960177572e-05\\
1.02522771981195	0.00012224968162311\\
1.10843914928567	1.52898537653101e-06\\
1.22569500209686	7.06217217611993e-05\\
1.18212541693221	1.31641350264999e-06\\
1.59832425596118	7.44899150806333e-07\\
2.10212996787622	7.70246330716873e-07\\
2.81311530429035	1.12608651292589e-07\\
3.99657268286855	5.32638018135856e-08\\
5.46397804739886	1.33608026864351e-08\\
7.06985608990456	4.16853485063484e-09\\
8.7929380516256	1.14338583045992e-09\\
12.4906584069235	1.302007390791e-10\\
16.0229192264476	3.18729043158328e-10\\
21.8545430038692	6.65245636355394e-12\\
29.2663641109433	3.00498737004773e-10\\
36.6605911364062	4.01348287937253e-10\\
51.6398816394326	2.17558415727126e-10\\
67.0670225568921	8.34357027912347e-10\\
89.3631052928142	1.12289200160376e-09\\
119.978461366658	1.147068218188e-09\\
};

\addplot [color=black, line width=0.5pt, solid]
  table[row sep=crcr]{%
0.27438139732593	69.143024572697\\
0.375778891364824	68.6930680566039\\
0.430757056740005	43.2983841758332\\
0.554608476223269	79.7825938878128\\
0.535169701279758	77.1434023842689\\
0.599850743445935	41.4190581387187\\
0.649170548427778	59.0316036272882\\
0.718909862790239	49.489979966796\\
0.774238955180106	47.9719347552226\\
0.825005996379058	48.3855438127671\\
0.877080383299957	40.9006432667609\\
0.942915071759922	3.45121625246258\\
1.09131457139624	1.13245703237228\\
1.01245325601131	3.22040654995179\\
1.10837130539746	2.52197411969148\\
1.16217961399229	2.21352813735734\\
1.23767395127908	0.806177018749089\\
1.65243095758646	0.114260779897472\\
2.17671831908061	0.000668852473520776\\
2.88591930119594	0.000541067913612636\\
3.75817436299542	3.91211231898225e-06\\
5.11025622475179	7.46575374011016e-06\\
6.69811090794635	4.05386460400337e-06\\
8.98501760789055	1.10795994601176e-06\\
12.0172881837261	3.45126860423761e-09\\
16.08929786679	5.04326047590098e-09\\
21.2792595530049	2.23367235818728e-09\\
27.3241956221479	1.86542381541699e-09\\
35.3904109568478	1.33511868227743e-10\\
46.8656922162947	1.02544417401873e-09\\
61.6319559122793	1.48270862254662e-09\\
81.8939858194264	1.51115919777567e-09\\
108.811054471136	1.3529577458371e-09\\
};

\addplot [color=colorclassyorange, solid, mark=10-pointed star,mark
size=1,mark options={solid},mark repeat=1]
  table[row sep=crcr]{%
2.57650974415349	10.0293607498953\\
3.13481901442043	6.20552505122641\\
3.7669244234095	9.25487900534435\\
4.26598736722786	9.29507189108535\\
4.96175525947736	8.31274025079756\\
5.48249192402565	11.5429362858658\\
6.28126875275616	14.645432835166\\
6.65028276792257	12.6869721388933\\
7.30898478225564	8.57584186661356\\
7.80098715859289	11.328605484234\\
8.65958128313206	5.16659098328052\\
8.98231766129509	10.4170433023926\\
9.69810100153789	1.95100496292573\\
10.3429325430621	1.95757567510364\\
10.8491560190287	4.35533925918163\\
11.3722918355884	4.35746803924623\\
12.3438334258432	0.503098588137472\\
16.2399749397885	0.000152399853704743\\
22.0767947785422	2.46411885198583e-05\\
29.8183987221333	3.59437191037593e-05\\
39.7271835150159	1.88838538228708e-06\\
53.5345200810404	7.08342984268029e-08\\
71.3219249532522	1.31552774895738e-07\\
94.9195785753804	5.87179385114212e-08\\
127.100798306424	1.25457511046534e-09\\
169.521438218384	3.17371906533026e-10\\
227.431741693702	2.98747693250334e-11\\
301.032387951766	3.11195069713222e-10\\
401.345652902863	4.06092048876872e-10\\
536.483354682144	2.16041406986278e-10\\
711.430332882337	8.3384410487497e-10\\
955.601760578921	1.12268150331829e-09\\
1311.75962550173	1.14693721187109e-09\\
};

\addplot [color=colorclassyorange, dashed, mark=triangle,mark size=1,mark
options={solid},mark repeat=1]
  table[row sep=crcr]{%
1.40777869203588	10.9664944266125\\
1.66284859675724	14.3550275155756\\
2.00975931329977	11.6107201964703\\
2.30861974611376	14.8852474187719\\
2.57905749229176	12.3961577836211\\
2.81711929447159	17.0629311837632\\
3.242743931049	12.7014888837083\\
3.55667411248337	21.2932110228841\\
3.81573029592844	18.1097839275809\\
4.19153650496506	20.9260509945046\\
4.51921017706354	17.1018659089637\\
4.74620021692032	18.7766130937429\\
5.1707285442953	17.1486838650711\\
5.41595538153632	16.3529525435747\\
5.7724140746649	17.395277791378\\
6.07634298627782	18.041755304413\\
6.29367637924974	21.359779347262\\
8.69734000925199	19.6541496426341\\
11.4954558103637	18.3614209023073\\
15.2387057187294	0.327511886318357\\
20.6126450200635	1.88838300418936e-06\\
27.0577900844056	7.08342984268029e-08\\
37.5416614007601	1.31552774451649e-07\\
49.9810907676176	5.87179407318672e-08\\
67.1037792347871	1.25457377819771e-09\\
89.5216897811194	3.17372794711446e-10\\
120.945039546082	2.98752134142433e-11\\
159.172730239192	3.11194181534802e-10\\
212.416599898475	4.06092492966081e-10\\
283.318356845043	2.16041406986278e-10\\
375.388963120242	8.33842772607341e-10\\
499.391357469148	1.12268061513987e-09\\
667.602455828676	1.1469376559603e-09\\
};

\addplot [color=colorclassyorange, densely dotted, mark=square,mark
size=1,mark options={solid},mark repeat=1]
  table[row sep=crcr]{%
0.657482925704638	1.10720383783545\\
0.767474078981693	16.9045006987546\\
0.92762017749282	15.601588622683\\
1.0648236344091	24.5583961411377\\
1.2170392026602	16.4003066629779\\
1.35188133203369	25.5672555682025\\
1.51018844090242	21.1069795121477\\
1.62281710038224	26.6706129707234\\
1.8561214249498	22.7901674174583\\
1.83359605328687	24.4796721974552\\
2.07858693789012	24.6195358045807\\
2.19534566911335	20.307999337318\\
2.34232675117534	24.4373986516784\\
2.4906475999324	23.2772585242079\\
2.52674535156793	20.2325084090028\\
2.78033019560354	29.0279988657505\\
2.93500075198646	29.387300246597\\
3.96794586417954	12.6165719901597\\
5.10519825274971	27.7657139669748\\
6.98929147265366	2.34689623148446\\
9.05648093751848	18.9626453575675\\
12.3231527475124	20.4786675537396\\
16.3145692948769	9.89931983490108\\
21.7580942111452	0.0207526901844393\\
28.4151781787638	6.68552546656542e-09\\
39.1944754061402	3.16660031529636e-10\\
52.3257863511906	2.98778779495024e-11\\
69.1760750029194	3.11196401980851e-10\\
92.7917159610557	4.06090716609242e-10\\
123.694412515456	2.16041406986278e-10\\
163.756366293355	8.3384366078576e-10\\
217.174728031562	1.1226819474075e-09\\
289.519253405027	1.1469376559603e-09\\
};

\end{axis}
\end{tikzpicture}%
\\
        \vspace{0.5cm}
        \hspace{-0.6cm}
%
\definecolor{mycolor1}{rgb}{0.00000,0.44706,0.74118}%
\definecolor{mycolor2}{rgb}{0.95000,0.32500,0.09800}%
\definecolor{mycolor3}{rgb}{1.00000,0.99000,0.97000}%
\begin{tikzpicture}
    \begin{axis}[%
    width=5.25in, height=0.75in,
    hide axis,
    xmin=0,
    xmax=60,
    ymin=0,
    ymax=0.4,
    legend style={draw=white!15!black,legend cell align=left, legend columns=-1, font=\scriptsize}
    ]
    \addlegendimage{colorclassyblue, line width=0.5pt, mark=o, mark size=1, mark options={solid}}
    \addlegendentry{IKS $\qquad$};
    \addlegendimage{black, line width=0.5pt, solid}
    \addlegendentry{BIKS $\qquad$};
    \addlegendimage{colorclassyorange, line width=0.5pt, mark=10-pointed star, mark size=1, mark options={solid}}
    \addlegendentry{AF50 $\qquad$};
    \addlegendimage{colorclassyorange, line width=0.5pt, dashed, mark=triangle, mark size=1, mark options={solid}}
    \addlegendentry{AF25 $\qquad$};
    \addlegendimage{colorclassyorange, line width=0.5pt, densely dotted, mark=square, mark size=1, mark options={solid}}
    \addlegendentry{AF10};
    \end{axis}
\end{tikzpicture}
    \caption{$\CC{L}^2$ error (top row) and error in the total energy (bottom row) for $\ve=10^{-2}$.
    First two columns present accuracy {\em vs} time step for $M=1000$ spatial grid points and
    $M=2000$ spatial grid points, respectively.
    The two subsequent columns illustrate the accuracy {\em vs} cost in seconds.}
    \label{fig:e1by100accuracy}
\end{figure}

\begin{figure}[tbh]
%
%
\definecolor{mycolor1}{rgb}{0.00000,0.44700,0.74100}%
\definecolor{mycolor2}{rgb}{0.85000,0.32500,0.09800}%
\definecolor{mycolor3}{rgb}{0.92900,0.69400,0.12500}%
\definecolor{mycolor4}{rgb}{0.49400,0.18400,0.55600}%
\begin{tikzpicture}

\begin{axis}[%
width=1.8in, height=1.4in, scale only axis, xmode=log, xmin=0.001, xmax=0.1,
xminorticks=true, ymode=log, ymin=1e-9, ymax=10, yminorticks=true, axis
background/.style={fill=white}, xticklabel style =
{xshift=-0.25ex,yshift=-0.1ex}, ylabel={$\CC{L}^2$ error}, xlabel={$\ve$
\mbox{}}, legend style={legend cell align=left, align=left,
draw=black,font=\scriptsize}, legend pos=south east, grid=both,
title={$\sigma=1, c=\Frac12$}]

\addplot [color=colorclassyblue,line width=0.5pt, solid, mark=o,mark
size=1,mark options={solid},mark repeat=1,solid]
  table[row sep=crcr]{%
0.1	1.61920256263195\\
0.0562341325190349	0.512769607532389\\
0.0316227766016838	0.124271342020701\\
0.0177827941003892	0.00923568997586365\\
0.01	0.000418925659363386\\
0.00562341325190349	5.89029538981187e-07\\
0.00316227766016838	1.36813826970419e-07\\
0.00177827941003892	2.21454138876002e-08\\
0.001	5.36536071613056e-09\\
}; \addlegendentry{IKS}

\addplot [color=black, line width=0.5pt, solid]
  table[row sep=crcr]{%
0.1	1.40098584909249\\
0.0562341325190349	1.92834185997349\\
0.0316227766016838	1.25300707147438\\
0.0177827941003892	0.520322729707463\\
0.01	0.836219649010207\\
0.00562341325190349	0.116915807602239\\
0.00316227766016838	0.00196097250528292\\
0.00177827941003892	0.000307518995576054\\
0.001	4.12411028451705e-06\\
}; \addlegendentry{BIKS}

\addplot [color=colorclassyorange, solid, mark=10-pointed star,mark
size=1,mark options={solid},mark repeat=1]
  table[row sep=crcr]{%
0.1	1.29126890639356\\
0.0562341325190349	1.99069863009746\\
0.0316227766016838	0.198383021468518\\
0.0177827941003892	0.456754943228247\\
0.01	0.332391236263864\\
0.00562341325190349	0.0120713218937939\\
0.00316227766016838	0.0908187107571974\\
0.00177827941003892	1.36064185492604\\
0.001	1.41421356238008\\
}; \addlegendentry{AF50}

\end{axis}
\end{tikzpicture}%
        \hspace{0.2in}
%
%
\definecolor{mycolor1}{rgb}{0.00000,0.44700,0.74100}%
\definecolor{mycolor2}{rgb}{0.85000,0.32500,0.09800}%
\definecolor{mycolor3}{rgb}{0.92900,0.69400,0.12500}%
\definecolor{mycolor4}{rgb}{0.49400,0.18400,0.55600}%
\begin{tikzpicture}

\begin{axis}[%
width=1.8in, height=1.4in, scale only axis, xmode=log, xmin=0.0005, xmax=0.1,
xminorticks=true, ymode=log, ymin=1e-5, ymax=10, yminorticks=true, axis
background/.style={fill=white}, xticklabel style =
{xshift=-0.25ex,yshift=-0.1ex}, ylabel={$\CC{L}^2$ error}, xlabel={$\ve$
\mbox{}}, legend style={legend cell align=left, align=left,
draw=black,font=\scriptsize}, legend pos=south east, grid=both,
title={$\sigma=\Frac12, c=2$}]

\addplot [color=colorclassyblue,line width=0.5pt, solid, mark=o,mark
size=1,mark options={solid},mark repeat=1,solid]
  table[row sep=crcr]{%
0.1	0.605131867802339\\
0.0562341325190349	0.673824423351313\\
0.0316227766016838	0.0701112374412092\\
0.0177827941003892	0.337041091924671\\
0.01	0.160638283734724\\
0.00562341325190349	0.00469487908586441\\
0.00316227766016838	0.00284895788179941\\
0.00177827941003892	0.0030641198741902\\
0.001	0.000322771706204646\\
0.000562341325190349	5.33339301379795e-05\\
}; \addlegendentry{IKS}

\addplot [color=black, line width=0.5pt, solid]
  table[row sep=crcr]{%
0.1	1.39568336491032\\
0.0562341325190349	1.40965024629419\\
0.0316227766016838	1.32819929673061\\
0.0177827941003892	1.33350498830158\\
0.01	1.29426017331692\\
0.00562341325190349	1.64530257396394\\
0.00316227766016838	1.2157713758019\\
0.00177827941003892	1.41753681822019\\
0.001	1.41673572311509\\
0.000562341325190349	1.45987454540412\\
}; \addlegendentry{BIKS}

\addplot [color=colorclassyorange, solid, mark=10-pointed star,mark
size=1,mark options={solid},mark repeat=1]
  table[row sep=crcr]{%
0.1	1.16430187463895\\
0.0562341325190349	1.56152787147618\\
0.0316227766016838	1.34534998603846\\
0.0177827941003892	1.67831359645337\\
0.01	1.41438538641521\\
0.00562341325190349	1.41416584638932\\
0.00316227766016838	1.4142121807698\\
0.00177827941003892	1.41473896673162\\
0.001	1.41421357888922\\
0.000562341325190349	1.41421356159475\\
}; \addlegendentry{AF50}

\end{axis}
\end{tikzpicture}%
    \caption{Global $\CC{L}^2$ errors under the time-step scaling $h = \ve^{\sigma}/c$ with $\sigma = 1$, $c=\Frac12$
    and $M = \lceil 20/\ve \rceil$ spatial grid points (left) and
    $\sigma = \Frac12$, $c=2$ and $M = \lceil 10/\ve \rceil$
    spatial grid points (right).}
    \label{fig:sigma}
\end{figure}
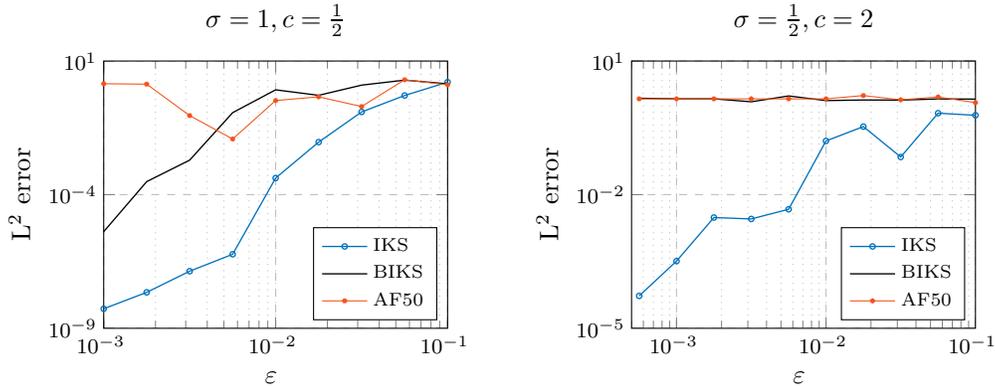

AF is accompanied by a postfix which refers to the number of Lanczos
iterations (10, 25 or 50) used in the method. In all cases presented here,
IKS and BIKS utilise two Lanczos iterations for the exponentiation of
$\WW^{[3]}$. Three Lanczos iterations are used for the exponentiation of
$W^{[2]}$ in experiments where time-steps tend to be small
(Figure~\ref{fig:e1by100accuracy}, Figure~\ref{fig:sigma} (left)) and five
iterations are used once large time-steps are involved
(Table~\ref{tab:largestep} and Figure~\ref{fig:sigma} (right)).

In all cases presented, the integrals in IKS \R{eq:Zass4} are discretised via
eleven Gauss--Legendre knots. In principle we can also use analytic or
asymptotic approximations, as well as highly oscillatory quadrature for the
integrals. Note that BIKS and AF use a fixed number of Gauss--Legendre knots
-- three and four, respectively.

Since the potentials are available in their analytic form, we use analytic
derivatives in our implementation. We also exploit the special form of the
potential to speed-up computations further without affecting accuracy. In
particular, for a potential of the form $V_{\mathrm{D}}(x) + f(t) x$, where
time dependence is scalar, the integrals $\momS{\star}{j}{k}$ and
$\LfS{f}{\star}{a}{b}$ can be reduced to a scalar form, the spatial
derivatives do not have be computed at each step and certain terms such as
$\dx^2 \mu_{o,1,1}(h), \dx^3 \mu_{o,1,1}(h), \dx^2 \mu_{e,2,1}(h)$ and $\dx^3
\mu_{o,3,1}(h)$ vanish. Relevant optimisation steps are applied to all methods
under consideration.


\begin{table}
\centering
\begin{tabular}{|c|c|c|c|c|c|c|c|}
  \hline
  $N$ &  & 30 & 40 & 50 & 60 & 75 & 100 \\
  $h$ &  & 0.083 & 0.063 & 0.050 & 0.042 & 0.033 & 0.025 \\
   \hline
   & IKS & 0.4967 & 0.3937 & 0.1606 & 0.0447 & 0.0021 & 0.0004 \\
  $\CC{L}^2$ & BIKS & 1.3855 & 1.3971 & 1.2943 & 0.7883 & 1.8290 & 0.3783   \\
  error & AF50 & 1.4148 & 1.4143 & 1.4144 & 1.4132 & 1.4020 & 0.3513  \\
  \hline
  Error & IKS & 0.1609 & 0.0012 & 0.0012  & 0.00046 & 0.000016 & 0.000001 \\
  in & BIKS & 4.0501 & 11.420 & 15.620 & 0.98735 & 0.332851 & 0.000536  \\
  energy& AF50 & 11.686 & 9.2533 & 10.183 & 7.36785 & 11.23268 & 0.003209    \\
  \hline
   & IKS & 0.44 & 0.56 & 0.74 & 0.77 & 0.90 & 1.25 \\
  Cost & BIKS & 0.42 & 0.53 & 0.64 & 0.74 & 0.92 & 1.19   \\
  (s) & AF50 & 2.63 & 3.47 & 4.33 & 5.20 & 6.39 & 8.45  \\
  \hline
\end{tabular}
\caption{
Juxtaposition of $\CC{L}^2$ error,
error in total energy (an observable) and cost in seconds of
analysed methods for relatively large
time steps, $h$, for the case $\ve=10^{-2}$. Here $N$ is the number of time steps ($h=T/N$).
The number of spatial grid points is chosen to be $M=1000$.} \label{tab:largestep}
\end{table}

{\bf Large time steps.} From Table~\ref{tab:largestep} we can see that BIKS
and AF50 do not provide physically meaningful solutions when the time step
$h$ is larger than $\ve$. This is due to the poor resolution of the highly
oscillatory potential. In the case of AF50, the large number of Lanczos
iterations required when using large-time steps is also a very important
factor in cost {\em vs} accuracy considerations. Due to the discretisation of the
integrals in IKS via eleven Gauss-Legendre knots and the separation of scales
in the asymptotic splitting \R{eq:Zass4}, it behaves exceptionally well for
much larger time steps, requiring very few Lanczos iterations for achieving
this accuracy.

{\bf Small time steps.} Having shown that IKS method works very well for
large time steps, in Figure~\ref{fig:e1by100accuracy} we consider the case of
smaller time steps, presenting two cases, with $M=1000$ and $M=2000$ spatial
grid points. The first two columns of Figure~\ref{fig:e1by100accuracy}
presents the accuracy. Once again, IKS seems to have significantly better
accuracy than BIKS (not specialised for highly oscillatory potentials) and
AF10, AF25 and AF50 (not specialised for semiclassical regime). The third and
fourth columns present the efficiency (accuracy {\em vs} cost) of methods. IKS is
not only more accurate, but far less costly.

Note that in the case of $M=1000$, although we achieve a high accuracy in
total energy ($\sim 10^{-9}$), the accuracy in $\CC{L}^2$ norm saturates at a
relatively low level ($\sim 10^{-6}$). This is inevitable since the space
resolution becomes inadequate beyond this stage. Analysing the same behaviour
 for $M=2000$ (columns 2 and 4 of Figure~\ref{fig:e1by100accuracy})
we observe accuracy of $\sim 10^{-9}$ both in terms of total energy and
$\CC{L}^2$ norm.

{\bf Semiclassical regime.} As the semiclassical parameter $\ve$ becomes
smaller, the effects discussed above become more pronounced. In particular,
choosing a time step as small or smaller than $\ve$ becomes infeasible. As we
seek larger time-steps, such as $h = \sqrt{\ve}/2$ used in
Figure~\ref{fig:sigma} (right), for instance, AF50 and BIKS are not able to
attain any meaningful accuracy, in contrast to the presented method.

{\bf Reference solutions.} The reference solutions for these experiments were
generated by using the sixth-order schemes BIKS and AF10 with very small time
steps. In the case of $\ve=10^{-2}$ experiments, we use $M=5000$ grid points
and $h=2.5 \times 10^{-6}$ as time step, which corresponds to $N=10^6$ time
steps. For Figure~\ref{fig:sigma}, where $\ve$ is varied, we use BIKS for
generating reference solutions, using $M=\lceil 20/\ve \rceil$ spatial grid
points.


\section{Conclusions}
To conclude, in this paper we have introduced a new numerical approach for
the approximation of \schr equation in semiclassical regime. This methodology
is especially powerful in the case of highly oscillatory, time-dependent
potentials. We have outlined a procedure for deriving arbitrarily high-order
numerical schemes, and have presented second, fourth and sixth order schemes
of this class of methods.

As demonstrated with the help of numerical examples, these schemes are more
effective in the semiclassical regime than
%
\begin{enumerate}
\item methods based on direct exponentiation via Lanczos iterations, which
    are costly in the semiclassical regime, and
\item the asymptotic Magnus--Zassenhaus splittings presented in
    \cite{BIKS}, which are specialised for the semiclassical regime, once
    the potential features rapid temporal oscillations.
\end{enumerate}
Let us briefly summarise the reason for the effectiveness of the current
method in the regime under consideration. In the semiclassical regime, a
decreasing $\ve$ results in increasing oscillations in both space and time,
consequently requiring an increase in the resolution of spatial and temporal
grids. Additionally, the highly oscillatory potential requires high temporal
grid resolution.
\begin{enumerate}
\item All methods compared in this paper require a very fine spatial grid
    in the semiclassical regime. The consequent imposition of a fine
    temporal grid in the case of exponential propagators, however, can be
    seen to arise from an interplay between the growth of the spectral
    radius (due to a fine spatial grid) and the need to keep Lanczos
    iterations small. In this case, resorting to an asymptotic splitting
    such as the symmetric Zassenhaus splitting of \cite{bader14eaf} and
    \cite{BIKS} allows us to separate the large scales dominating the
    spectral radius, whereby the only terms requiring exponentiation via
    Lanczos iterations are asymptotically small (in increasing powers of
    $\ve$) and impose considerably more relaxed constraints on the time step.
\item By preserving integrals till the very last stage, the new method
    decouples the task of approximating the integrals of the potential from
    the task of approximating the time evolution of the solution. Thus, it
    allows us to resort to higher accuracy or specialised quadrature
    methods for the approximation of the integrals of the potential,
    $V(t,x)$, without affecting the time step of the scheme. This is
    particularly effective in applications where the time-dependent part of
    the potential is oscillatory but not excessively expensive to evaluate,
    even if the time-independent part is. In all other known methods
    (including \cite{BIKS}), highly oscillatory potentials require a very
    high resolution in time since the sampling of the potential (via
    quadrature methods or otherwise) is fixed at the outset and linked to
    the
    time-propagation. 
\end{enumerate}
%



It goes without saying that, like in virtually every mathematical paper, our
exposition is incomplete. In particular, no attempt has been made in this
paper to optimise the performance of the Magnus--Zassenhaus approach or to
exploit serendipitous connections among the terms {\em \'a la\/}
\cite{blanes00iho} to reduce the cost further. This will be matter for future
research.

\section*{Acknowledgments}
The work of Karolina Kropielnicka in this project was financed by The National Center for
Science, based on grant no. 2016/22/M/ST1/00257.

\bibliographystyle{agsm}
\bibliography{MZSemiclassical}

\appendix

\section{The algebra of anti-commutators}
\label{sec:algebra}
The vector field of the \schr equation \R{eq:schr} is a linear combination of
the action of two operators, $\dx^2$ and the multiplication by the
interaction potential $V(t)$, for any $t \geq 0$. Since the Magnus expansion
requires nested commutation, the focus of our interest is the Lie algebra
generated by $\dx^2$ and $V(\cdot)$,
\begin{displaymath}
  \GG{F}=\CC{LA} \{\dx^2, V(\cdot) \},
\end{displaymath}
i.e.\ the  linear-space closure of all nested commutators of $\dx^2$ and
$V(\cdot)$.

{\bf Simplifying commutators}. To simplify commutators we could study their
action on functions. For example, using the chain rule we find,
\begin{displaymath}
  [\dx^2,V]u=\dx^2(Vu) - V(\dx^2u)=(\dx^2V)u+2(\dx V)\dx u,
\end{displaymath}
which implies that $[\dx^2,V]=(\dx^2V)+2(\dx V)\dx$. 
Straightforward discretisation of these simplified commutators,
\[ (\dx^2V)+2(\dx V)\dx \leadsto \DDD_{(\dx^2V)} + 2 \DDD_{(\dx V)} \K_1,\]
are not skew-Hermitian (or even Hermitian, for that matter). In practice,
this results in a loss of unitarity and stability when such matrices are
exponentiated.

In the methods presented here we circumvent the problem by working with {\it
anti-commutators}, that is the differential operators of the form (\ref{eq:angbra}),
\begin{displaymath}
  \ang{f}{k} := \Frac12 \left( f \circ \dx^k + \dx^k \circ f \right), \quad  k \geq 0,\ f \in \CC{C}_p^{\infty}(I;\BB{R}),
\end{displaymath}
which are inherently symmetrised. The action of this differential operator on
$u$, for example, is
$$ \ang{f}{k} u = \Frac12 \left( f \dx^k u + \dx^k (f u) \right) $$
and the discretisation of this operator, given by (\ref{eq:angbradisc}), is
\begin{displaymath}
 \ang{f}{k} \leadsto \Frac12 \left( \DDD_f \K_k + \K_k \DDD_f \right).
\end{displaymath}
It is  simple matter to verify that the discretisation of $\ang{f}{k}$
is skew-Hermitian for odd $k$ and Hermitian for even $k$. This is the reason
why the choice of algebra of anti-commutators $\ang{\cdot}{k}$ seems to be
optimal for our purposes.

Moreover, the commutators of these symmetrised differential operators can be
simplified using the following rules,
\begin{Eqnarray}
\label{eq:id}\left[ \Ang{1}{f}, \Ang{0}{g} \right] & =& \Ang{0}{f (\dx g)},\\
\nonumber\left[ \Ang{1}{f}, \Ang{1}{g} \right] & =& \Ang{1}{f (\dx g)-(\dx f) g},\\
\nonumber\left[ \Ang{2}{f}, \Ang{0}{g} \right] & =& 2 \Ang{1}{f (\dx g)},\\
\nonumber\left[ \Ang{2}{f}, \Ang{1}{g} \right] & =& \Ang{2}{2 f (\dx g) -
(\dx f) g} - \Frac12 \Ang{0}{2(\dx f) (\dx^2 g) + f (\dx^3
g)},\\
\nonumber\left[ \Ang{2}{f}, \Ang{2}{g} \right] & =& 2 \Ang{3}{f (\dx g) -
(\dx f) g} + \Ang{1}{2(\dx^2 f) (\dx g) - 2(\dx f) (\dx^2 g) + (\dx^3 f) g -
f (\dx^3 g)},\\
\nonumber\left[ \Ang{3}{f}, \Ang{0}{g} \right] & =& 3 \Ang{2}{f (\dx g)} - \Frac12 \Ang{0}{3(\dx f) (\dx^2 g) + f (\dx^3 g)},\\
\nonumber\left[ \Ang{4}{f}, \Ang{0}{g} \right] & =& 4 \Ang{3}{f (\dx g)} -2
\Ang{1}{3(\dx f) (\dx^2 g) + f (\dx^3 g)}.
\end{Eqnarray}

There is rich algebraic theory underlying these differential operators
\cite{psalgebra}, including a general formula for  \R{eq:id}.  In principle,
however, the above rules can be verified by application of the chain rule.

\begin{remark}
Note that, by definition \R{eq:angbra},
\begin{enumerate}
    \item these brackets are linear, so that $\Ang{2}{2 f (\dx g) - (\dx f)
        g} = 2 \Ang{2}{f (\dx g)} - \Ang{2}{(\dx f) g}$,
    \item $\ang{f}{0} = f$; and
    \item $\ang{1}{2} = \dx^2$.
\end{enumerate}
\end{remark}

With this new notation in place and using \R{eq:id}, we can now simplify
commutators using the terminology of anti-commutators,
    \begin{Eqnarray}
        \nonumber[\ii\dx^2, \ii V ] & =& -[ \ang{1}{2}, \ang{V}{0} ] = - 2 \ang{\dx V}{1},\\
        \nonumber[\ii V,[\ii\dx^2,\ii V]] & = & -\ii [ \ang{V}{0},[ \ang{1}{2},\ang{V}{0}]] = 2\ii \ang{(\dx V)^2}{0},\\
        \nonumber[\ii \dx^2,[\ii \dx^2,\ii V]] & = & -\ii [\ang{1}{2},[\ang{1}{2}, \ang{V}{0} ]] = \ii \ang{\dx^4 V}{0} - 4 \ii \ang{\dx^2 V}{2}.
    \end{Eqnarray}
Straightforward discretisations of these operators preserve the symmetries
that are crucial for retaining unitarity.

\section{Symmetric Zassenhaus splitting}
\label{sec:Zassenhaus}
The symmetric Zassenhaus splitting algorithm has crucial advantages over
conventional splittings \cite{bader14eaf}. It provides neat separation of
terms with differing scales and structures, each of which is easy to
exponentiate separately either on account of the structure or of the size. It
achieves this separation by a recursive application of the symmetric {\em
Baker--Campbell--Hausdorff formula} (usually known in an abbreviated form as
the sBCH formula),
\begin{equation}
  \label{eq:sBCH}
  \ee^{\frac12 X}\ee^{Y}\ee^{\frac12 X}=\ee^{\sBCH(X,Y)},
\end{equation}
for $X$ and $Y$ in a Lie algebra $\GG{g}$, where
\begin{equation}
\label{eq:sBCHcoeffs}
\CC{sBCH}(\dt X,\dt Y) = \dt (X+Y)-\dt^3(\Frac{1}{24}[[Y,X],X]+\Frac{1}{12}[[Y,X],Y]) +\O{\dt^5}.
\end{equation}
The expansion \R{eq:sBCHcoeffs} features terms of the Hall basis
\cite{reutenauer93fla}, such as $[[Y,X],X]$ and $[[Y,X],Y]$, and can be
computed to an arbitrary power of $\dt$ using an algorithm from
\cite{casas09eac}. Coefficients and terms of the Hall basis for sufficiently
high degree sBCH expansions are also available in a tabular form
\cite{murua10tsb}. Since \R{eq:sBCH} is palindromic, only odd powers of $\dt$
feature in the expansion \R{eq:sBCHcoeffs}.

For instance, consider the task of approximating $\ee^{\mathcal{W}^{[0]}}$,
where $\mathcal{W}^{[0]}=X+Y$, to an accuracy of $\O{h^5}$ under the
assumption $X,Y = \O{\dt}$. Using the sBCH formula \R{eq:sBCH}, and with the
choice of $W^{[0]}=X$ as our first exponent, we write
\begin{equation}\label{eq:2:sbch_inverted}
  \ee^{\mathcal{W}^{[0]}} =\ee^{\frac12 W^{[0]}} \ee^{\sBCH(- W^{[0]}, \mathcal{W}^{[0]})} \ee^{\frac12  W^{[0]}},
\end{equation}
\begin{Eqnarray}
\nonumber \mathcal{W}^{[1]} := \CC{sBCH}(-X,X+Y) &=& (-X + (X+Y)) - \Frac{1}{24}[[X+Y,X],X]\\
\nonumber &&+\Frac{1}{12}[[X+Y,X],X+Y] + \O{h^5} \\
\label{eq:sBCHshifted}&=&Y + \Frac{1}{24}[[Y,X],X] + \Frac{1}{12}[[Y,X],Y] +
\O{h^5}.
\end{Eqnarray}
Note that $\mathcal{W}^{[1]} = Y + \O{\dt^3}$ and thus we have extracted
$W^{[0]}=X$ from the exponent $\WW^{[0]} = X + Y$ at the cost of correction
terms in form of higher-order commutators $Z = \Frac{1}{24}[[Y,X],X] +
\Frac{1}{12}[[Y,X],Y] = \O{h^3}$ appearing\footnote{Since we are expanding to
an accuracy of $\O{h^5}$, we can safely discard here all $\O{h^5}$ terms.} in
$\WW^{[1]}$. At the next stage,
 we extract the largest term from $\WW^{[1]}=Y+Z$, which is
$W^{[1]}=Y=\O{h}$,
\begin{equation}\label{eq:2:sbch_inverted2}
  \ee^{\mathcal{W}^{[1]}} =\ee^{\frac12 W^{[1]}} \ee^{\sBCH(- W^{[1]}, \mathcal{W}^{[1]})} \ee^{\frac12  W^{[1]}},
\end{equation}
and $\mathcal{W}^{[2]} = \sBCH(- W^{[1]}, \mathcal{W}^{[1]}) = Z + \O{h^5}$
features grade three and higher commutators of $Y$ and $Z$ (or grade five and
higher commutators of $X$ and $Y$) as corrections, which can be discarded.
Thus, by combining \R{eq:2:sbch_inverted} and \R{eq:2:sbch_inverted2}, the
recursive algorithm terminates with the splitting
\begin{align*}
	 \exp(X+Y) = \ee^{\frac12 W^{[0]}} \ee^{\frac12 W^{[1]}} \ee^{\WW^{[2]}}
									 \ee^{\frac12 W^{[1]}} \ee^{\frac12 W^{[0]}} + \O{h^5},
\end{align*}
where $W^{[0]}$ and $W^{[1]}$ are $\O{h}$ and $\WW^{[2]}$ is $\O{h^3}$. A
sixth-order splitting would feature an additional exponent,
$\WW^{[3]}=\O{h^5}$. These exponential splittings are, thus, characterised by
an expansion in increasing powers of the parameter $h$ and we say that they
are {\em asymptotic splittings} in $h$.

We emphasize that, in principle, we are free to choose the elements $W^{[k]}$
 that we want to extract in each step of the algorithm,
\begin{equation}\label{eq:2:sbch_invertedK}
  \ee^{\mathcal{W}^{[k]}} =\ee^{\frac12 W^{[k]}} \ee^{\mathcal{W}^{[k+1]}} \ee^{\frac12  W^{[k]}}, \quad \mathcal{W}^{[k+1]} = \sBCH(- W^{[k]}, \mathcal{W}^{[k]}).
\end{equation}
  This choice can afford a great deal of
flexibility -- it could be based on some structural property which allows for
trivial exponentiation of $W^{[k]}$ when extracted separately, a small
spectral radius which makes the term amenable to effective exponentiation
through Krylov subspace methods, a combination of both, or some other
criteria. As long as the terms are decreasing in size, the convergence of the
procedure is guaranteed. This can lead to many variants of the splitting,
some of which could prove to have more favourable properties than others.

Analysing the sizes of terms in powers of $h$ might seem natural in some
cases such as ODEs but is not the only choice. In the case of the \schr
equation in the semiclassical regime, for instance, the semiclassical
parameter $\ve$ is a more natural choice due to the considerations outlined
in Subsection~2.2. In this case the exponents in a symmetric Zassenhaus
splitting,
\begin{align*}
	 \exp(X+Y) = \ee^{\frac12 W^{[0]}} \ee^{\frac12 W^{[1]}}
									\cdots \ee^{\frac12 W^{[s]}} \ee^{\WW^{[s+1]}} \ee^{\frac12 W^{[s]}}
									\cdots \ee^{\frac12 W^{[1]}} \ee^{\frac12 W^{[0]}} + \O{\ve^r},
\end{align*}
are expanded in powers of $\ve$.

\end{document}